# EFFECTIVE CHOICE AND BOUNDEDNESS PRINCIPLES IN COMPUTABLE ANALYSIS

VASCO BRATTKA AND GUIDO GHERARDI

**Abstract.** In this paper we study a new approach to classify mathematical theorems according to their computational content. Basically, we are asking the question which theorems can be continuously or computably transferred into each other? For this purpose theorems are considered via their realizers which are operations with certain input and output data. The technical tool to express continuous or computable relations between such operations is Weihrauch reducibility and the partially ordered degree structure induced by it. We have identified certain choice principles such as co-finite choice, discrete choice, interval choice, compact choice and closed choice, which are cornerstones among Weihrauch degrees and it turns out that certain core theorems in analysis can be classified naturally in this structure. In particular, we study theorems such as the Intermediate Value Theorem, the Baire Category Theorem, the Banach Inverse Mapping Theorem, the Closed Graph Theorem and the Uniform Boundedness Theorem. We also explore how existing classifications of the Hahn-Banach Theorem and Weak Kőnig's Lemma fit into this picture. Well-known omniscience principles from constructive mathematics such as LPO and LLPO can also naturally be considered as Weihrauch degrees and they play an important role in our classification. Based on this we compare the results of our classification with existing classifications in constructive and reverse mathematics and we claim that in a certain sense our classification is finer and sheds some new light on the computational content of the respective theorems. Our classification scheme does not require any particular logical framework or axiomatic setting, but it can be carried out in the framework of classical mathematics using tools of topology, computability theory and computable analysis. We develop a number of separation techniques based on a new parallelization principle, on certain invariance properties of Weihrauch reducibility, on the Low Basis Theorem of Jockusch and Soare and based on the Baire Category Theorem. Finally, we present a number of metatheorems that allow to derive upper bounds for the classification of the Weihrauch degree of many theorems and we discuss the Brouwer Fixed Point Theorem as an example.

**§1. Introduction.** The purpose of this paper is to propose a new approach to classify mathematical theorems according to their computational content and according to their logical complexity.

2000 *Mathematics Subject Classification.* 03F60,03D30,03B30,03E15.

*Key words and phrases.* Computable analysis, constructive analysis, reverse mathematics, effective descriptive set theory.

This project has been supported by the Italian Ministero degli Affari Esteri and the National Research Foundation of South Africa (NRF). An extended abstract version of this paper appeared in [11].

October 1, 2018

### 1.1. Realizability of theorems and Weihrauch reducibility.

The basic idea is to interpret theorems, which are typically $\Pi_2$–theorems of the form

$$(\forall x \in X)(\exists y \in Y)(x, y) \in A,$$

as operations $F :\subseteq X \rightrightarrows Y, x \mapsto \{y \in Y : (x, y) \in A\}$ that map certain input data $X$ into certain output data $Y$. In other words, we are representing theorems by their realizers or multi-valued Skolem functions, which is a very natural approach for many typical theorems. For instance, the Intermediate Value Theorem states that

$$(\forall f \in \mathcal{C}[0, 1], f(0) \cdot f(1) < 0)(\exists x \in [0, 1])f(x) = 0$$

and hence it is natural to consider the partial multi-valued operation

$$\mathsf{IVT} :\subseteq \mathcal{C}[0, 1] \rightrightarrows [0, 1], f \mapsto \{x \in [0, 1] : f(x) = 0\}$$

with $\mathrm{dom}(\mathsf{IVT}) := \{f \in \mathcal{C}[0, 1] : f(0) \cdot f(1) < 0\}$ as a representative of this theorem. It follows from the Intermediate Value Theorem itself that this operation is well-defined. The goal of our study is to understand the computational content of theorems like the Intermediate Value Theorem and to analyze how they compare to other theorems. In order to understand the relation of two theorems $T$ and $T'$ to each other we will ask the question whether a realizer $G$ of $T'$ can be computably or continuously transformed into a realizer $F$ of $T$. In other words, we consider theorems as points in a space (represented by their realizers) and we study whether these points can be computably or continuously transferred into each other. This study is carried out entirely in the domain of classical logic and using tools from topology, computability theory and computable analysis [52].

In fact the technical tool to express the relation of realizers to each other is a reducibility that Weihrauch introduced in the 1990s in two unpublished papers [50, 51] and which since then has been studied by several others (see for instance [23, 2, 5, 36, 20, 12, 39]). Basically, the idea is to say that a single-valued function $F$ is Weihrauch reducible to $G$, in symbols $F \leq_{\mathrm{W}} G$, if there are computable functions $H$ and $K$ such that

$$F = H\langle \mathrm{id}, GK \rangle.$$

Here $K$ can be considered as an input adaption and $H$ as an output adaption. The output adaption has direct access to the input, since in many cases the input cannot be looped through $G$. Here and in the following $\langle\ \rangle$ denotes suitable finite or infinite tupling functions. This reducibility can be extended to sets of functions and to multi-valued functions on represented spaces. The resulting structure has been studied in [12] and among other things it has been proved that parallelization is a closure operator for Weihrauch reducibility. To parallelize a multi-valued function $F$ just means to consider

$$\widehat{F}(p_0, p_1, p_2, \ldots) := F(p_0) \times F(p_1) \times F(p_2) \times \ldots,$$

i.e. to take countably many instances of $F$ in parallel. If $f$ is defined on Baire space $\mathbb{N}^{\mathbb{N}}$, then we sometimes compose parallelization with an infinite tupling function. This is convenient, but does not affect the operation in any essential way. The resulting parallelized partial order forms a lattice into which Turing and Medvedev degrees can be embedded. In Section 2 we will summarize the definition of Weihrauch reducibility and some relevant results. In this paper we will mainly study the non-parallelized Weihrauch degrees of theorems since they allow a finer classification of computational properties. Nevertheless, the closure operator of parallelization will play an important role.

**1.2. Effective choice and boundedness principles.** A characterization of the Weihrauch degree of theorems is typically achieved by showing that the degree is identical to the degree of some other known principle. We have identified certain choice principles that turned out to be crucial cornerstones in our classification. We can formulate these choice principles as follows:

- $(C_F)$ Any non-empty co-finite set $A \subseteq \mathbb{N}$ has a member $x \in A$.
- $(C_{\mathbb{N}})$ Any non-empty set $A \subseteq \mathbb{N}$ has a member $x \in A$.
- $(C_I^-)$ Any proper closed interval $I \subseteq [0, 1]$ has a member $x \in I$.
- $(C_I)$ Any non-empty closed interval $I \subseteq [0, 1]$ has a member $x \in I$.
- $(C_K)$ Any non-empty compact set $K \subseteq [0, 1]$ has a member $x \in K$.
- $(C_A)$ Any non-empty closed sets $A \subseteq \mathbb{R}$ has a member $x \in A$.
- $(C)$ Any set $A \subseteq \mathbb{N}$ has a characteristic function $\mathrm{cf}_A : \mathbb{N} \to \{0, 1\}$.

By a proper closed interval, we mean an interval that has more than one point. We will refer to these principles as *co-finite choice*, *discrete choice*, *proper interval choice*, *interval choice*, *compact choice*, *closed choice* and *countable choice*. As they stand, these principles are trivially correct in classical mathematics. They only become interesting as soon as one considers them from an algorithmic point of view. Given a set $A$ by negative information, typically by some form of enumeration of its complement, how difficult is it to actually find a member of $A$ (or the characteristic function $\mathrm{cf}_A$ of $A$)? It is not too hard to see that from this point of view these principles are algorithmically unsolvable and they form some hierarchy of principles of different degree of unsolvability. If we consider the realizers of these choice principles in the above mentioned sense, then they all correspond to discontinuous operations of different degree of discontinuity.

Often it is more convenient to consider these choice principles as boundedness principles and in particular the principles of interval choice have equivalent boundedness versions. The benefit of the boundedness principles is that the negative information about the represented set $A$ is given explicitly in form of one or two bounds and indeed many problems in analysis can be reduced to finding such bounds. In particular, we will consider the following boundedness principles:

- $(B_F)$ For any sequence $(q_n)_{n \in \mathbb{N}}$ of rational numbers bounded from above, there exists a real number $x \in \mathbb{R}$ with $\sup_{n \in \mathbb{N}} q_n \leq x$.
- $(B_I^-)$ For any two sequences $(q_n)_{n \in \mathbb{N}}$ and $(r_n)_{n \in \mathbb{N}}$ of rational numbers such that $\sup_{n \in \mathbb{N}} q_n < \inf_{n \in \mathbb{N}} r_n$, there exists a real number $x \in \mathbb{R}$ such that $\sup_{n \in \mathbb{N}} q_n \leq x \leq \inf_{n \in \mathbb{N}} r_n$.
- $(B_I)$ For any two sequences $(q_n)_{n \in \mathbb{N}}$ and $(r_n)_{n \in \mathbb{N}}$ of rational numbers such that $\sup_{n \in \mathbb{N}} q_n \leq \inf_{n \in \mathbb{N}} r_n$, there exists a real number $x \in \mathbb{R}$ such that $\sup_{n \in \mathbb{N}} q_n \leq x \leq \inf_{n \in \mathbb{N}} r_n$.
- $(B_I^+)$ For any two sequences $(q_n)_{n \in \mathbb{N}}$ and $(r_n)_{n \in \mathbb{N}}$ of rational numbers such that $\sup_{n \in \mathbb{N}} q_n \leq \inf_{n \in \mathbb{N}} r_n$, there exists a real number $x \in \mathbb{R}$ such that $\sup_{n \in \mathbb{N}} q_n \leq x \leq \inf_{n \in \mathbb{N}} r_n$. We allow the case $r_n = \infty$.

- (B) For any sequence $(q_n)_{n \in \mathbb{N}}$ of rational numbers bounded from above, there exists a real number $x \in \mathbb{R}$ with $\sup_{n \in \mathbb{N}} q_n = x$.

Once again, these statements are classically correct and even trivial as they stand. However, finding a bound algorithmically (given the guarantee that it exists) is a different story and once again these principles represent operations of different degree of discontinuity. In Section 3 we will prove the equivalence of certain choice and boundedness principles and we will compare them to omniscience principles. Omniscience principles have been introduced by Brouwer and Bishop [1, 15] as non-acceptable principles in the intuitionistic framework of constructive analysis. Intuitionistic reasoning does neither allow the law of the excluded middle ($A \vee \neg A$) nor de Morgan's law $\neg(A \wedge B) \iff (\neg A \vee \neg B)$. If these laws are applied to simple existential statements in first-order arithmetic, i.e. statements of the form $A = (\exists n \in \mathbb{N})P(n)$, then the law of excluded middle and de Morgan's law translate into LPO and LLPO, respectively. The abbreviations stand for *limited principle of omniscience* and *lesser limited principle of omniscience*, respectively and the principles are typically formulated as follows:

- (LPO) For any sequence $p \in \mathbb{N}^{\mathbb{N}}$ there exists an $n \in \mathbb{N}$ such that $p(n) = 0$ or $p(n) \neq 0$ for all $n \in \mathbb{N}$.
- (LLPO) For any sequence $p \in \mathbb{N}^{\mathbb{N}}$ such that $p(k) \neq 0$ for at most one $k \in \mathbb{N}$, it follows $p(2n) = 0$ for all $n \in \mathbb{N}$ or $p(2n + 1) = 0$ for all $n \in \mathbb{N}$.

Once again the realizers of these statements correspond to discontinuous operations of different degree of discontinuity [51]. Other principles that are refuted in constructive mathematics, such as Markov's principle have continuous realizers and hence they are not problematic from our point of view.

The parallelizations $\widehat{\mathsf{LPO}}$ and $\widehat{\mathsf{LLPO}}$ turned out to be particularly important cornerstones in our classification scheme, since $\widehat{\mathsf{LPO}}$ is a $\mathbf{\Sigma}_2^0$–complete operation in the effective Borel hierarchy [5], i.e. it is complete among all limit computable operations with respect to Weihrauch reducibility and similarly $\widehat{\mathsf{LLPO}}$ is complete among all weakly computable operations [20, 12]. Limit computable operations are exactly the effectively $\mathbf{\Sigma}_2^0$–measurable operations and these are exactly those operations that can be computed on a Turing machine that is allowed to revise its output. We have defined weakly computable operations exactly by the above mentioned completeness property in [12]. In Section 3 we will also show how the choice and boundedness principles are related to the omniscience principles and their parallelizations.

In Section 4 we develop a number of separation techniques to show that certain principles are strictly stronger than others or even incomparable. These techniques include:

- The *Parallelization Principle* that states that nothing above LPO can be reducible to anything below $\widehat{\mathsf{LLPO}}$. This principle excludes several reductions between choice principles.
- The *Mind Change Principle* that considers the number of mind changes required to compute the corresponding function on a limit machine. This principle helps to separate decision problems from others.

- The *Computable Invariance Principle* that considers realizers that map computable values to computable values. This principle helps to separate problems with non-uniform computable solutions from problems that do not admit such solutions. It can be combined with well-known counterexamples from computable analysis, such as Specker's non-computable but left computable real number or Kreisel and Lacombe's co-c.e. compact set $A \subseteq [0, 1]$ that contains no computable points.
- The *Low Invariance Principle* that considers realizers that map computable values to low values. This principle can be combined, for instance, with the Low Basis Theorem of Jockusch and Soare.
- The *Baire Category Principle* that can help to separate non-discrete operations from discrete ones.

We apply these separation techniques to separate all choice principles (except co-finite and discrete choice, which are equivalent) and omniscience principles from each other. Figure 1 illustrates the relation between the choice principles and other results discussed in this paper.

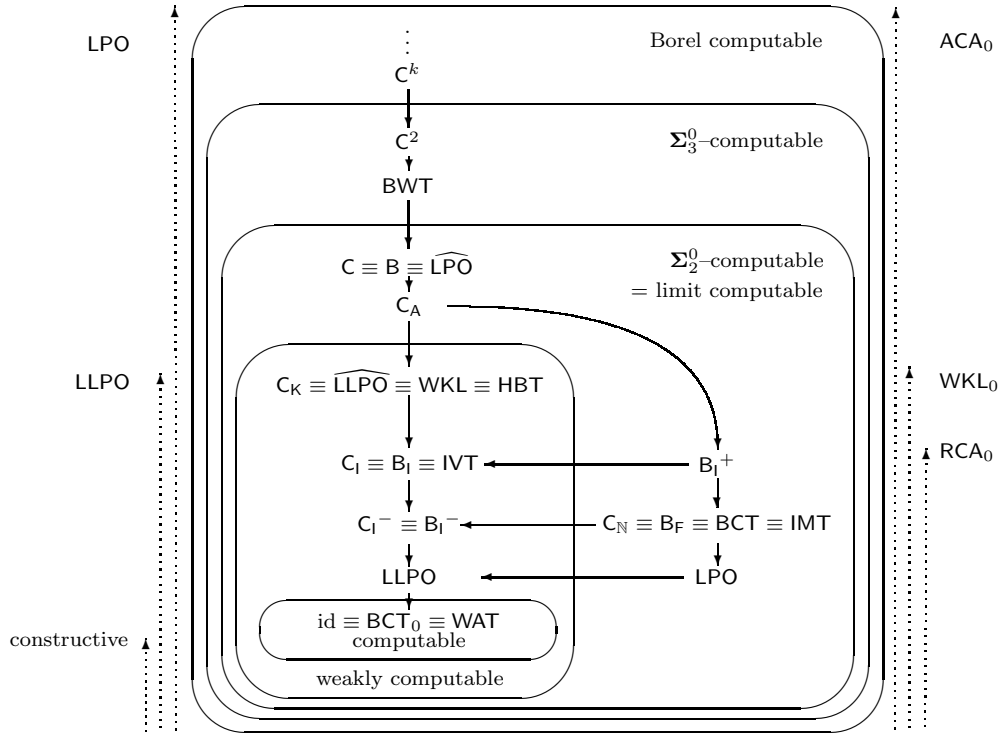

FIGURE 1. Constructive, computable and reverse mathematics

**1.3. Theorems in functional analysis.** As a case study we analyze a number of theorems from analysis and functional analysis and we classify their

Weihrauch degree. In particular, we will consider in Sections 5, 6, 7 and 8 the following theorems:

- **(BCT$_0$)** Given a sequence $(A_i)_{i \in \mathbb{N}}$ of closed nowhere dense subsets of a complete separable metric space $X$, there exists a point $x \in X \setminus \bigcup_{i \in \mathbb{N}} A_i$ (Baire Category Theorem).
- **(BCT)** Given a sequence $(A_i)_{i \in \mathbb{N}}$ of closed subsets of a complete separable metric space $X$ with $X = \bigcup_{i=0}^{\infty} A_i$, there is some $n \in \mathbb{N}$ such that $A_n$ is somewhere dense (Baire Category Theorem).
- **(IMT)** Any bijective linear bounded operator $T : X \to Y$ on separable Banach spaces $X$ and $Y$ has a bounded inverse $T^{-1} : Y \to X$ (Banach Inverse Mapping Theorem).
- **(OMT)** Any surjective linear bounded operator $T : X \to Y$ on separable Banach spaces $X$ and $Y$ is open, i.e. $T(U)$ is open for any open $U \subseteq X$ (Open Mapping Theorem).
- **(CGT)** Any linear operator $T : X \to Y$ with a closed graph$(T) \subseteq X \times Y$ is bounded (Closed Graph Theorem).
- **(UBT)** Any sequence $(T_i)_{i \in \mathbb{N}}$ of linear bounded operators that is pointwise bounded, i.e. such that $\sup\{\|T_i x\| : i \in \mathbb{N}\}$ exists for all $x \in X$, is uniformly bounded, i.e. $\sup\{\|T_i\| : i \in \mathbb{N}\}$ exists (Uniform Boundedness Theorem).
- **(HBT)** Any bounded linear functional $f : Y \to \mathbb{R}$, defined on some closed subspace $Y$ of a Banach space $X$ has a bounded linear extension $g : X \to \mathbb{R}$ with the same norm $\|g\| = \|f\|$ (Hahn-Banach Theorem).
- **(IVT)** For any continuous function $f : [0,1] \to \mathbb{R}$ with $f(0) \cdot f(1) < 0$ there exists a $x \in [0,1]$ with $f(x) = 0$ (Intermediate Value Theorem).
- **(BFT)** Any continuous function $f : [0,1]^n \to [0,1]^n$ has a fixed point $x \in [0,1]^n$, i.e. $f(x) = x$ (Brouwer Fixed Point Theorem).
- **(BWT)** Any sequence $(x_i)_{i \in \mathbb{N}}$ of numbers in $[0,1]^n$ has a convergent subsequence (Bolzano-Weierstraß Theorem).
- **(WAT)** For any continuous function $f : [0,1] \to \mathbb{R}$ and any $n \in \mathbb{N}$ there exists a rational polynomial $p \in \mathbb{Q}[x]$ such that $\|f - p\| = \sup_{x \in [0,1]} |f(x) - p(x)| < 2^{-n}$ (Weierstraß Approximation Theorem).
- **(WKL)** Any infinite binary tree has an infinite path (Weak Kőnig's Lemma).

The Baire Category Theorem is an example of a theorem for which it matters which version is realized. In the formulation $\mathsf{BCT}_0$ it leads to a continuous and even computable realizer, whereas the version $\mathsf{BCT}$ is discontinuous. The realizers of the given theorems are operations of different degree of discontinuity and our aim is classify the computational Weihrauch degree of these results. The benefit of such a classification is that practically all purely computability theoretic questions of interest about a theorem in computable analysis can be answered by such a classification. Typical questions are:

1. Is the theorem uniformly computable, i.e. can we compute the output information $y \in Y$ uniformly from the input information $x \in X$?
2. Is the theorem non-uniformly computable, i.e. does there exist a computable output information $y \in Y$ for any computable input information $x \in X$?
3. If there is no uniform solution, is there a uniform computation of a certain effective Borel complexity?

4. If there is no non-uniform computable solution, is there always a non-uniform result of a certain arithmetical complexity or Turing degree?

Answers to questions of this type can be derived from the classification of the Weihrauch degree of a theorem. In the diagram of Figure 1 we summarize some of our results. The arrows in the diagram are pointing into the direction of computations and implicit logical implications and hence in the inverse direction of the corresponding reductions. No arrow in the diagram can be inverted and no arrows can be added (except those that follow by transitivity).

In Section 8 we provide a number of metatheorems that allow to determine upper bounds of the Weihrauch degree of many theorems straightforwardly, just because of the mere topological form of the statement. For instance, any classical result of the form

$$(\forall x \in X)(\exists y \in Y)(x, y) \in A$$

with a co-c.e. closed $A \subseteq X \times Y$ and a co-c.e. compact $Y$ has a realizer that is reducible to compact choice $\mathsf{C_K}$. The table in Figure 2 summarizes the topological types of metatheorems and the corresponding version of computability. We illustrate that these metatheorems are useful and we show that one directly gets upper bounds for theorems such as the Brouwer Fixed Point Theorem and the Peano Existence Theorem for the initial value problem of ordinary differential equations.

| metatheorem | computability | unique case |
| --- | --- | --- |
| open | computable | computable |
| compact | weakly computable | computable |
| locally compact | limit computable | non-uniformly computable |

FIGURE 2. Types of metatheorems, choice and computability

**1.4. Related approaches.** Several other approaches to classify mathematical theorems from a logical point of view have already been presented and studied intensively. In our context, the most relevant approaches are constructive analysis as studied by Bishop, Bridges and Ishihara [1, 15] and many others and reverse mathematics as proposed by Friedman and Simpson [48]. In computable analysis theorems have been classified according to the Borel complexity by the authors of this papers and others [5, 19, 10, 20]:

- **Constructive mathematics:** in constructive mathematics theorems have been proved to be intuitionistically equivalent to certain principles which are not acceptable from the constructive point of view. Such principles are, for instance, the Limited Principle of Omniscience $\mathsf{LPO}$ and the Lesser Limited Principle of Omniscience $\mathsf{LLPO}$.
- **Reverse mathematics:** in reverse mathematics theorems have been classified according to which comprehension axioms are required to prove these theorems in second-order arithmetic. Such comprehension axioms are, for instance, the Recursive Comprehension Axiom $\mathsf{RCA_0}$ and the Arithmetic Comprehension Axiom $\mathsf{ACA_0}$. Another system of wide importance is the system of Weak Kőnig's Lemma $\mathsf{WKL_0}$.

- **Computable mathematics:** in computable mathematics theorems have been classified according to their Borel complexity, i.e. according to the complexity that realizers of those theorems have in terms of effective Borel measurability.

There are several further related approaches and variants of the aforementioned ones that we cannot list here completely. We just mention intuitionistic reverse mathematics as studied by Ishihara [26, 27] and others and uniform reverse mathematics as proposed by Kohlenbach [31, 30]. It is not too difficult to recognize that all these various approaches mentioned here are related in some form or the other, although they are expressed in terms that appear to be different on the first sight and they produce different and sometimes incompatible results.

We claim that our approach, which is a refinement of the third approach, sheds some new light on the computational status of theorems and perhaps on all the above approaches. In some respects, our classification scheme is finer than those proposed in other approaches and it distinguishes certain aspects that have not been captured yet. The arrows at the side in Figure 1 give some rough indication where results are located in other approaches. In the Conclusion we will add some further comments on how our approach compares to other existing ones.

**§2. Weihrauch reducibility, omniscience principles and weak computability.** In this section we briefly recall some definitions and results from [12] on Weihrauch reducibility and omniscience principles as far as they are relevant for the present study. More details and many further results can be found in the aforementioned source. We assume that the reader has some basic familiarity with concepts from computable analysis and otherwise we refer the reader for all undefined concepts to [52]. In a first step we define Weihrauch reducibility for sets of functions on Baire space, as it was already considered by Weihrauch [50, 51].

DEFINITION 2.1 (Weihrauch reducibility). Let $\mathcal{F}$ and $\mathcal{G}$ be sets of functions of type $f :\subseteq \mathbb{N}^{\mathbb{N}} \to \mathbb{N}^{\mathbb{N}}$. We say that $\mathcal{F}$ is *Weihrauch reducible* to $\mathcal{G}$, in symbols $\mathcal{F} \leq_{\mathrm{W}} \mathcal{G}$, if there are computable functions $H, K :\subseteq \mathbb{N}^{\mathbb{N}} \to \mathbb{N}^{\mathbb{N}}$ such that

$$(\forall G \in \mathcal{G})(\exists F \in \mathcal{F}) \ F = H\langle \mathrm{id}, GK \rangle.$$

Analogously, we define $\mathcal{F} \leq_{\mathrm{sW}} \mathcal{G}$ using the equation $F = HGK$ and in this case we say that $\mathcal{F}$ is *strongly Weihrauch reducible* to $\mathcal{G}$.

Here $\langle \ \rangle : \mathbb{N}^{\mathbb{N}} \times \mathbb{N}^{\mathbb{N}} \to \mathbb{N}^{\mathbb{N}}$ denotes a computable standard pairing function [52]. This reducibility is derived from Weihrauch reducibility of single functions (which we did not formalize here) in the same way as Medvedev reducibility is derived from Turing reducibility in classical computability theory [43]. We denote the induced equivalence relations by $\equiv_{\mathrm{W}}$ and $\equiv_{\mathrm{sW}}$, respectively.

In the next step we define the concept of a realizer of a multi-valued function as it is used in computable analysis [52]. We recall that a *representation* $\delta_X :\subseteq \mathbb{N}^{\mathbb{N}} \to X$ of a set $X$ is a surjective (and potentially partial) map. In general, the inclusion symbol "$\subseteq$" indicates partiality in this paper. In this situation we say that $(X, \delta_X)$ is a *represented space*.

DEFINITION 2.2 (Realizer). Let $(X, \delta_X)$ and $(Y, \delta_Y)$ be represented spaces and let $f :\subseteq X \rightrightarrows Y$ be a multi-valued function. Then $F :\subseteq \mathbb{N}^{\mathbb{N}} \to \mathbb{N}^{\mathbb{N}}$ is called *realizer* of $f$ with respect to $(\delta_X, \delta_Y)$, in symbols $F \vdash f$, if

$$\delta_Y F(p) \in f\delta_X(p)$$

for all $p \in \operatorname{dom}(f\delta_X)$.

Usually, we do not mention the representations explicitly since they will be clear from the context. A multi-valued function $f :\subseteq X \rightrightarrows Y$ on represented spaces is called *continuous* or *computable*, if it has a continuous or computable realizer, respectively. Using reducibility for sets and the concept of a realizer we can now define Weihrauch reducibility for multi-valued functions.

DEFINITION 2.3 (Realizer reducibility). Let $f$ and $g$ be multi-valued functions on represented spaces. Then $f$ is said to be *Weihrauch reducible* to $g$, in symbols $f \leq_{\mathrm{W}} g$, if and only if $\{F : F \vdash f\} \leq_{\mathrm{W}} \{G : G \vdash g\}$. Analogously, we define $f \leq_{\mathrm{sW}} g$ with the help of $\leq_{\mathrm{sW}}$ on sets.

That is, $f \leq_{\mathrm{W}} g$ holds if any realizer of $g$ computes some realizer of $f$ with some fixed uniform translations $H, K$. This reducibility has already been used in [5] for single-valued maps and in [20] for multi-valued maps. We mention that we also write $f <_{\mathrm{W}} g$ if and only if $f \leq_{\mathrm{W}} g$ and $g \not\leq_{\mathrm{W}} f$. Moreover, we write $f \mid_{\mathrm{W}} g$ if $f \not\leq_{\mathrm{W}} g$ and $g \not\leq_{\mathrm{W}} f$. Analogous notation is used for $\leq_{\mathrm{sW}}$. It is clear that Weihrauch reducibility and its strong version form preorders, i.e. both relations are reflexive and transitive. It is also clear that strong reducibility is actually stronger than weak reducibility in the sense that $f \leq_{\mathrm{sW}} g$ implies $f \leq_{\mathrm{W}} g$ and both reducibilities preserve computability and continuity from right to the left. By a *Weihrauch degree* we mean the equivalence class of all multi-valued functions on represented spaces with respect to Weihrauch reducibility, where we only consider functions with at least one computable point in the domain and defined on a certain fixed set of underlying represented spaces. Among all Weihrauch degrees there is a least one, namely the degree of the computable multi-valued functions. Weihrauch reducibility and strong Weihrauch reducibility are both invariant under equivalent representations, i.e. if we replace representations by computably equivalent ones, this does not affect the reducibility relations.

One can show that the product of multi-valued functions $f \times g$ and the direct sum $f \oplus g$ are both monotone operations with respect to strong and ordinary Weihrauch reducibility and hence both operations can be extended to Weihrauch degrees. This turns the structure of partially ordered Weihrauch degrees into a lower-semi lattice with the direct sum operation as greatest lower bound operation. It turns out that a very important operation on this lower semi-lattice is parallelization, which can be understood as countably infinite product operation.

DEFINITION 2.4 (Parallelization). Let $f :\subseteq X \rightrightarrows Y$ be a multi-valued function. Then we define the *parallelization* $\widehat{f} :\subseteq X^{\mathbb{N}} \rightrightarrows Y^{\mathbb{N}}$ of $f$ by

$$\widehat{f}(x_i)_{i \in \mathbb{N}} := \underset{i=0}{\overset{\infty}{\times}} f(x_i)$$

for all $(x_i)_{i \in \mathbb{N}} \in X^{\mathbb{N}}$.

We mention that parallelization acts as a closure operator with respect to Weihrauch reducibility.

PROPOSITION 2.5 (Parallelization). *Let $f$ and $g$ be multi-valued functions on represented spaces. Then*

1. $f \leq_{\mathrm{W}} \widehat{f}$                                          *(extensive)*
2. $f \leq_{\mathrm{W}} g \Longrightarrow \widehat{f} \leq_{\mathrm{W}} \widehat{g}$                              *(increasing)*
3. $\widehat{f} \equiv_{\mathrm{W}} \widehat{\widehat{f}}$                                          *(idempotent)*

*An analogous result holds for strong Weihrauch reducibility.*

Hence, we can define in a natural way a parallelized version of Weihrauch reducibility and of Weihrauch degrees, just by replacing equivalence classes by their closures. In some sense, the relation between ordinary Weihrauch reducibility and its parallelized version is similar to the relation between many-one reducibility and Turing reducibility in classical computability theory. The parallelized Weihrauch degrees even form a lattice with the product operations $f \times g$ as least upper bound operation and the direct sum as $f \oplus g$ as greatest lower bound operation. We can attach a virtual greatest element to this lattice that acts like a multi-valued function without realizer. One can prove that Medvedev degrees can be embedded into this lattice such that least upper bounds and greatest lower bounds are preserved and as a special case one obtains an embedding of Turing degrees.

Very helpful cornerstones in the Weihrauch semi-lattice are the omniscience principles that we already mentioned in the introduction. Formally, we consider them as maps.

DEFINITION 2.6 (Omniscience principles). We define:

- $\mathsf{LPO} : \mathbb{N}^{\mathbb{N}} \to \mathbb{N}, \quad \mathsf{LPO}(p) = \begin{cases} 0 & \text{if } (\exists n \in \mathbb{N}) \ p(n) = 0 \\ 1 & \text{otherwise} \end{cases}$,

- $\mathsf{LLPO} :\subseteq \mathbb{N}^{\mathbb{N}} \rightrightarrows \mathbb{N}, \mathsf{LLPO}(p) \ni \begin{cases} 0 & \text{if } (\forall n \in \mathbb{N}) \ p(2n) = 0 \\ 1 & \text{if } (\forall n \in \mathbb{N}) \ p(2n+1) = 0 \end{cases}$,

where $\mathrm{dom}(\mathsf{LLPO}) := \{p \in \mathbb{N}^{\mathbb{N}} : p(k) \neq 0 \text{ for at most one } k\}$.

One should notice that the definition of $\mathsf{LLPO}$ implies that $\mathsf{LLPO}(0^{\mathbb{N}}) = \{0, 1\}$. The following result summarizes the relations between the omniscience principles and their parallelizations to each other.

THEOREM 2.7 (Omniscience principles). *We obtain*

$$\mathsf{LLPO} <_{\mathrm{W}} \mathsf{LPO} \mid_{\mathrm{W}} \widehat{\mathsf{LLPO}} <_{\mathrm{W}} \widehat{\mathsf{LPO}}.$$

The parallelization $\mathsf{C} := \widehat{\mathsf{LPO}}$ of $\mathsf{LPO}$ is known to be $\boldsymbol{\Sigma}_2^0$–complete among all effectively $\boldsymbol{\Sigma}_2^0$–measurable maps (see [5]). The effectively $\boldsymbol{\Sigma}_2^0$–measurable maps are also called limit computable. Similarly, $\widetilde{\mathsf{LLPO}}$ plays a significant role and it can be proved that the class of all multi-valued functions that are reducible to $\widetilde{\mathsf{LLPO}}$ are closed under composition. An equivalent result has already been obtained in [20], expressed in terms of Weak König's Lemma. We do not mention the straightforward formalization of $\mathsf{WKL}$ as multi-valued function here, but we mention the result from [12].

PROPOSITION 2.8. $\mathsf{WKL} \equiv_{\mathrm{W}} \widehat{\mathsf{LLPO}}$.

From now on we will use $\mathsf{WKL}$ and $\widehat{\mathsf{LLPO}}$ interchangeably, depending on which one fits better into the context. We take the nice composition property of functions below $\mathsf{WKL} \equiv_{\mathrm{W}} \widehat{\mathsf{LLPO}}$ as a reason to call these functions *weakly computable*. This class of functions is closely related to compact choice, which we can define in the following generalized sense.

DEFINITION 2.9 (Compact choice). Let $X$ be a computable metric space. The multi-valued operation

$$\mathsf{C}_{\mathcal{K}(X)} :\subseteq \mathcal{K}_-(X) \rightrightarrows X, A \mapsto A$$

with $\mathrm{dom}(\mathsf{C}_{\mathcal{K}(X)}) := \{A \subseteq X : A \neq \emptyset \text{ compact}\}$ is called *compact choice* of $X$.

Here $\mathcal{K}_-(X)$ denotes the set of compact subsets of $X$, represented by enumerations of finite rational open covers (which are not necessarily minimal) (see [13]). It can be proved that $\widehat{\mathsf{LLPO}}$ is equivalent to compact choice for a large class of computable metric spaces, which has essentially been proved in [20].

THEOREM 2.10 (Compact choice). *Let $X$ be a computable metric space. Then $\mathsf{C}_{\mathcal{K}(X)} \leq_{\mathrm{sW}} \widehat{\mathsf{LLPO}}$. If there is a computable embedding $\iota : \{0,1\}^{\mathbb{N}} \hookrightarrow X$, then $\mathsf{C}_{\mathcal{K}(X)} \equiv_{\mathrm{sW}} \widehat{\mathsf{LLPO}}$.*

One can also prove that on computable metric spaces a multi-valued function is weakly computable if and only if it has an upper semi-computable compact-valued selector. This implies the following important result, which states that single-valued weakly computable functions are automatically fully computable.

COROLLARY 2.11 (Weak computability). *Let $X$ be a represented space and $Y$ a computable metric space. Any weakly computable single-valued operation $f :\subseteq X \to Y$ is computable.*

**§3. Choice and boundedness principles.** In this section we study choice principles and boundedness principles. Both types of principles are closely related to each other and they are also related to the omniscience principles mentioned earlier. In some sense most of the boundedness principles are just variants of the choice principles that are more convenient for some applications.

By $\mathcal{A}(X)$ or $\mathcal{A}_-(X)$ we denote the set of closed subsets of a metric space $X$. The index "$-$" indicates that we assume that the hyperspace $\mathcal{A}_-(X)$ is equipped with the lower Fell topology and a corresponding negative information representation $\psi_-$ explained below. All choice principles are restrictions of the multi-valued choice map

$$\text{Choice} :\subseteq \mathcal{A}_-(X) \rightrightarrows X, A \mapsto A,$$

which is defined for non-empty closed sets $A \subseteq X$ and maps any such set in a multi-valued way to the set of its members. That is, the input is a non-empty closed set $A \in \mathcal{A}_-(X)$ and the output is one of the (possibly many) points $x \in A$. We can define restrictions of the choice map by specifying the respective domains and ranges.

**Definition 3.1** (Choice principles). We define multi-valued operations as restrictions of the respective choice maps as follows:

1. $\mathsf{C_F} :\subseteq \mathcal{A}_-(\mathbb{N}) \rightrightarrows \mathbb{N}$, $\mathrm{dom}(\mathsf{C_F}) := \{A \subseteq \mathbb{N} : A \text{ co-finite}\}$.
2. $\mathsf{C_{\mathbb{N}}} :\subseteq \mathcal{A}_-(\mathbb{N}) \rightrightarrows \mathbb{N}$, $\mathrm{dom}(\mathsf{C_{\mathbb{N}}}) := \{A \subseteq \mathbb{N} : A \neq \emptyset\}$.
3. $\mathsf{C_I} :\subseteq \mathcal{A}_-[0,1] \rightrightarrows [0,1]$, $\mathrm{dom}(\mathsf{C_I}) := \{[a,b] : 0 \leq a \leq b \leq 1\}$.
4. $\mathsf{C_I}^- :\subseteq \mathcal{A}_-[0,1] \rightrightarrows [0,1]$, $\mathrm{dom}(\mathsf{C_I}^-) := \{[a,b] : 0 \leq a < b \leq 1\}$.
5. $\mathsf{C_K} :\subseteq \mathcal{A}_-([0,1]) \rightrightarrows [0,1]$, $\mathrm{dom}(\mathsf{C_K}) := \{K \subseteq [0,1] : K \neq \emptyset \text{ compact}\}$.
6. $\mathsf{C_A} :\subseteq \mathcal{A}_-(\mathbb{R}) \rightrightarrows \mathbb{R}$, $\mathrm{dom}(\mathsf{C_A}) := \{A \subseteq \mathbb{R} : A \neq \emptyset \text{ closed}\}$.

We refer to these operations as *co-finite choice*, *discrete choice*, *interval choice*, *proper interval choice*, *compact choice* and *closed choice*, respectively.

Whenever $X$ is a computable metric space, we will use the representation $\psi_-$ of $\mathcal{A}_-(X)$ that represents a closed set $A \subseteq X$ by negative information. There are different equivalent ways of characterizing this representation (see [13] for details). The intuition is that a name $p$ of a closed set $A$ is an enumeration of rational open balls $B(x_i, r_i)$ that exhaust the complement of $A$, i.e. $X \setminus A = \bigcup_{i=0}^{\infty} B(x_i, r_i)$. Here a rational open ball is a ball with center in the dense set and a rational radius. The $\psi_-$–computable members in $\mathcal{A}_-(X)$ are called *co-c.e. closed* sets.

In the case of $\mathcal{A}_-(\mathbb{N})$ we can consider $\psi_-$ just as an enumeration of the complement, i.e. $\psi_-(p) = A$ if and only if

$$\mathbb{N} \setminus A = \{n \in \mathbb{N} : n + 1 \in \mathrm{range}(p)\} =: \mathrm{range}(p) - 1.$$

Here $\mathrm{range}(p) - 1$ is used opposed to $\mathrm{range}(p)$ in order to leave the sequence $0^{\mathbb{N}}$ as name for $\mathbb{N}$.

For practical purposes it is often more convenient to handle these choice principles in form of the closely related boundedness principles that we define now.

**Definition 3.2** (Boundedness principles). We define the following multi-valued operations (with their maximal domains, if not stated otherwise):

1. $\mathsf{B_F} : \mathbb{R}_< \rightrightarrows \mathbb{R}, x \mapsto [x, \infty)$.
2. $\mathsf{B_I} :\subseteq \mathbb{R}_< \times \mathbb{R}_> \rightrightarrows \mathbb{R}, (x, y) \mapsto [x, y]$, $\mathrm{dom}(\mathsf{B_I}) := \{(x, y) : x \leq y\}$.
3. $\mathsf{B_I}^- :\subseteq \mathbb{R}_< \times \mathbb{R}_> \rightrightarrows \mathbb{R}, (x, y) \mapsto [x, y]$, $\mathrm{dom}(\mathsf{B_I}^-) := \{(x, y) : x < y\}$.
4. $\mathsf{B_I}^+ :\subseteq \mathbb{R}_< \times \overline{\mathbb{R}_>} \rightrightarrows \mathbb{R}, (x, y) \mapsto [x, y]$, $\mathrm{dom}(\mathsf{B_I}^+) := \{(x, y) : x \leq y\}$.
5. $\mathsf{B} : \mathbb{R}_< \rightrightarrows \mathbb{R}, x \mapsto x$.

Here we assume that $\mathbb{R}, \mathbb{R}_<$ and $\mathbb{R}_>$ are equipped with the ordinary Cauchy representation $\rho$ of the reals and the left and right representations $\rho_<$ and $\rho_>$, respectively (see [52] for details). A name $p \in \mathbb{N}^{\mathbb{N}}$ of a real number with respect to the representation $\rho_<$ is a sequence of rational numbers $(q_n)_{n \in \mathbb{N}}$ such that $x = \sup_{n \in \mathbb{N}} q_n$. Without loss of generality one can assume that the sequence is strictly increasing, whenever this is helpful. The representation $\rho_>$ is defined analogously using the infimum instead of the supremum and without loss of generality we can assume that the corresponding sequence is strictly decreasing, whenever required. The representation $\rho = \rho_< \sqcap \rho_>$ is the join of $\rho_<$ and $\rho_>$, which means that it contains both types of information simultaneously. Here $\overline{\mathbb{R}_>} = \mathbb{R} \cup \{\infty\}$ and it is represented by the extension of $\rho_>$ that represents a point by a decreasing sequence of numbers in $\mathbb{Q} \cup \{\infty\}$. The sequence that has constant value $\infty$ is the only name of $\infty$. We mention that $[a, \infty]$ is meant to be the set $\{x \in \mathbb{R} : a \le x < \infty\}$ in case of $\mathsf{B_I}^+$ since the result is required to be an ordinary real number. The boundedness principle $\mathsf{B_I}^-$ is just the restriction of $\mathsf{B_I}$ to the non-degenerate case and $\mathsf{B_I}^+$ is the extension of $\mathsf{B_I}$ to the case where the right-hand side bound is allowed to be $\infty$. Sometimes it is convenient to assume that $\mathrm{range}(\mathsf{B_F}) = \mathbb{N}$ and this is possible without loss of generality since the operation $U : \mathbb{R} \rightrightarrows \mathbb{N}, x \mapsto \{n \in \mathbb{N} : x \le n\}$ is computable (see [52]).

It is clear the the choice and the boundedness principles are closely related to each other and as a first result we prove that co-finite choice $\mathsf{C_F}$ is strongly equivalent to $\mathsf{B_F}$. Perhaps it is surprising that both operations are equivalent to discrete choice $\mathsf{C_\mathbb{N}}$, however not in the strong sense since this result uses the direct access to the input in an essential way.

PROPOSITION 3.3 (Discrete choice). $\mathsf{B_F} \equiv_{\mathrm{sW}} \mathsf{C_F} \equiv_{\mathrm{W}} \mathsf{C_\mathbb{N}}$.

PROOF. We first show $\mathsf{B_F} \le_{\mathrm{sW}} \mathsf{C_F}$. Given a bounded sequence $(q_i)_{i \in \mathbb{N}}$ of rational numbers we want to compute a real number $x$ with $\sup_{i \in \mathbb{N}} q_i \le x$. Now we consider the following co-finite set of natural numbers

$$A = \mathbb{N} \setminus \{n \in \mathbb{N} : (\exists i)\, n \le q_i\}.$$

Any point $k \in A$ satisfies $k \ge \sup_{i \in \mathbb{N}} q_i$. Given $(q_i)_{i \in \mathbb{N}}$ we can generate a $\psi_{--}$ name of $A$ and hence with the help of some realizer of $\mathsf{C_F}$, we can find a point $k \in A$, which is a desired upper bound. Thus $\mathsf{B_F} \le_{\mathrm{sW}} \mathsf{C_F}$.

We now prove $\mathsf{C_F} \le_{\mathrm{sW}} \mathsf{B_F}$. Given a sequence $p \in \mathbb{N}^{\mathbb{N}}$ such that the set $\mathrm{range}(p) - 1 = \mathbb{N} \setminus A$ is finite, we can consider $p$ as an enumeration of rational numbers and since $\mathbb{N} \setminus A$ is finite it follows that $\max \mathrm{range}(p)$ exists. We can use any realizer of $\mathsf{B_F}$ to find a natural number $n \ge \max \mathrm{range}(p)$ and any such $n$ is in $A$. Thus $\mathsf{C_F} \le_{\mathrm{sW}} \mathsf{B_F}$.

It is clear that $\mathsf{C_F} \le_{\mathrm{sW}} \mathsf{C_\mathbb{N}}$, because the former is a restriction of the latter. We prove $\mathsf{C_\mathbb{N}} \le_{\mathrm{W}} \mathsf{B_F}$. Let $p \in \mathbb{N}^{\mathbb{N}}$ be a sequence such that $\mathrm{range}(p) - 1 = \mathbb{N} \setminus A$ and $A$ is non-empty. The goal is to find a number $k \in A$. For this purpose we inspect the sequence $p$ and we try to find a suitable number (whose successor is) not enumerated in $p$. Starting with the candidate $c = 0$, we inspect $p$ in stages $n = 0, 1, 2, ...$ and whenever the current candidate $c$ is enumerated, i.e. whenever $p(n) = c + 1$, then we replace the current candidate by the smallest number which has not yet been enumerated (i.e. $c = \min(\mathbb{N} \setminus \{p(0) - 1, ..., p(n) - 1\})$. At the same time we record the stages in which we have to change our candidate in a sequence $r$, i.e. $r(n) := n$ in stages $n$ where a new candidate $c$ has been chosen and $r(n) := r(n-1)$ in all other stages $n > 1$. The sequence $r$ is increasing and eventually constant, since there is a smallest missing number in $\mathrm{range}(p) - 1 = \mathbb{N} \setminus A$. In fact, this smallest number is also the smallest natural number, which is not among the numbers

$$p(0) - 1, p(1) - 1, ..., p(\max \mathrm{range}(r)) - 1.$$

The function $F : \mathbb{N}^{\mathbb{N}} \to \mathbb{N}^{\mathbb{N}}$ that maps each $p$ to the corresponding $r = F(p)$ is computable and with the help of a realizer of $\mathsf{B_F}$ we can find a natural number $m \ge \max \mathrm{range}(r)$. With direct access to the original input $p$ it is then possible to inspect the numbers $p(0), ..., p(m)$ in order to find the smallest number $k \in \mathbb{N}$ that is not among the numbers $p(0) - 1, ..., p(m) - 1$. This number is a suitable output since $k \in A$. Altogether, this shows $\mathsf{C}_\mathbb{N} \leq_W \mathsf{B}_\mathsf{F}$. ⊣

We mention that discrete choice $\mathsf{C}_\mathbb{N}$ has another interesting property: it is idempotent. This means that $\mathsf{C}_\mathbb{N} \times \mathsf{C}_\mathbb{N} \equiv_W \mathsf{C}_\mathbb{N}$ or, in other words, one can use one instance of $\mathsf{C}_\mathbb{N}$ to solve two requests simultaneously. The relevance and properties of idempotency are discussed in [12].

PROPOSITION 3.4 (Idempotency of discrete choice). $\mathsf{C}_\mathbb{N} \equiv_W \mathsf{C}_\mathbb{N} \times \mathsf{C}_\mathbb{N}$.

PROOF. Since $\max : \mathbb{R}_< \times \mathbb{R}_< \to \mathbb{R}_<$ is computable and any $z \in \mathsf{B}_\mathsf{F}(\max(x, y))$ satisfies $z \in \mathsf{B}_\mathsf{F}(x) \cap \mathsf{B}_\mathsf{F}(y)$ we obtain $\mathsf{B}_\mathsf{F} \equiv_W \mathsf{B}_\mathsf{F} \times \mathsf{B}_\mathsf{F}$. By the previous Proposition 3.3 the same holds true for $\mathsf{C}_\mathbb{N}$. ⊣

Sometimes it is useful to assume that $\mathsf{B}_\mathsf{I}$ or $\mathsf{B}_\mathsf{I}^+$ are restricted to the open interval $(0, 1)$ and this is possible by the following lemma.

LEMMA 3.5. *The multi-valued function*

$$\mathsf{B}_\mathsf{I}' :\subseteq (0, 1)_< \times (0, 1)_> \rightrightarrows (0, 1), (x, y) \mapsto \{z : x \leq z \leq y\}.$$

*is strongly equivalent to* $\mathsf{B}_\mathsf{I}$, *i.e.* $\mathsf{B}_\mathsf{I} \equiv_{sW} \mathsf{B}_\mathsf{I}'$. *Analogously,* $\mathsf{B}_\mathsf{I}^- \equiv_{sW} \mathsf{B}_\mathsf{I}'^-$ *for the restriction* $\mathsf{B}_\mathsf{I}'^-$ *of* $\mathsf{B}_\mathsf{I}'$ *to* $\mathrm{dom}(\mathsf{B}_\mathsf{I}'^-) := \{(x, y) : x < y\}$.

PROOF. One can use the strictly monotone computable and bijective function

$$f : (0, 1) \to \mathbb{R}, x \mapsto \tan\left(\pi x - \frac{\pi}{2}\right)$$

in order to identify $(0, 1)$ with $\mathbb{R}$. The inverse $f^{-1}$ of $f$ is computable as well. We obtain $\mathsf{B}_\mathsf{I}'(x, y) = f^{-1}(\mathsf{B}_\mathsf{I}(f(x), f(y)))$ for all $x, y \in (0, 1)$ and $\mathsf{B}_\mathsf{I}(x, y) = f(\mathsf{B}_\mathsf{I}'(f^{-1}(x), f^{-1}(y)))$ for all $x, y \in \mathbb{R}$. Since $f$ is strictly monotone increasing, $f$ and its inverse $f^{-1}$ are also $(\rho_<, \rho_<)$– and $(\rho_>, \rho_>)$–computable. This implies $\mathsf{B}_\mathsf{I} \equiv_{sW} \mathsf{B}_\mathsf{I}'$. Analogously, one obtains $\mathsf{B}_\mathsf{I}^- \equiv_{sW} \mathsf{B}_\mathsf{I}'^-$. ⊣

Thus, wherever it is convenient, we can assume that the range of $\mathsf{B}_\mathsf{I}$ and $\mathsf{B}_\mathsf{I}^-$ is $(0, 1)$. The second main observation in this section is that interval choice $\mathsf{C}_\mathsf{I}$ is strongly equivalent to $\mathsf{B}_\mathsf{I}$ and correspondingly proper interval choice $\mathsf{C}_\mathsf{I}^-$ is strongly equivalent to $\mathsf{B}_\mathsf{I}^-$. This means that in case of interval choice the negative information can be provided in forms of left and right hand bounds (as in $\mathsf{B}_\mathsf{I}$ and $\mathsf{B}_\mathsf{I}^-$) or in form of intervals that exhaust the complement without explicit information on which side the negative information lies (as in $\mathsf{C}_\mathsf{I}$ and $\mathsf{C}_\mathsf{I}^-$). We note that in case of $\mathsf{B}_\mathsf{I}^+$ the negative information has to be given in form of two separately presented bounds because the underlying space is not compact.

PROPOSITION 3.6 (Interval choice). $\mathsf{B}_\mathsf{I} \equiv_{sW} \mathsf{C}_\mathsf{I}$, $\mathsf{B}_\mathsf{I}^- \equiv_{sW} \mathsf{C}_\mathsf{I}^-$, $\mathsf{B}_\mathsf{I}^+ \leq_{sW} \mathsf{C}_\mathsf{A}$.

PROOF. We prove $\mathsf{C}_\mathsf{I} \leq_{sW} \mathsf{B}_\mathsf{I}$. Let us assume that a non-empty interval $I = [a, b] \subseteq [0, 1]$ is given in the form $[0, 1] \setminus I = \bigcup_{i=0}^\infty J_i$, where any $J_i$ is an open interval in $[0, 1]$ given by its rational endpoints, possibly containing the points $0$ or $1$. We use this sequence $(J_i)_{i \in \mathbb{N}}$ to generate two sequences of real numbers

$(q_i)_{i \in \mathbb{N}}$ and $(r_i)_{i \in \mathbb{N}}$. We choose

$$q_n := \begin{cases} \sup(L) & \text{if } L = \{x \in \mathbb{R} : [0,x] \subseteq \bigcup_{i=0}^{n} J_i\} \neq \emptyset \\ 0 & \text{otherwise} \end{cases},$$

$$r_n := \begin{cases} \inf(U) & \text{if } U = \{y \in \mathbb{R} : [y,1] \subseteq \bigcup_{i=0}^{n} J_i\} \neq \emptyset \\ 1 & \text{otherwise} \end{cases}.$$

Then $\sup_{n \in \mathbb{N}} q_n = a$ and $\inf_{n \in \mathbb{N}} r_n = b$ because $[0,1] \setminus (a,b)$ is compact. Thus, given a realizer of $\mathsf{B_I}$, we can compute a real number $y \in \mathsf{B_I}(a,b)$, i.e. $a \leq y \leq b$. This shows $\mathsf{C_I} \leq_{\mathrm{sW}} \mathsf{B_I}$.

For the other direction we consider the boundedness principle $\mathsf{B_I}'$, which is $\mathsf{B_I}$ restricted to the open interval $(0,1)$, which is strongly equivalent to $\mathsf{B_I}$ according to Lemma 3.5. We prove $\mathsf{B_I}' \leq_{\mathrm{sW}} \mathsf{C_I}$. Given sequences $(q_i)_{i \in \mathbb{N}}$ and $(r_i)_{i \in \mathbb{N}}$ of rational numbers such that $a := \sup_{i \in \mathbb{N}} q_i \leq \inf_{i \in \mathbb{N}} r_i =: b$ with $a, b \in (0,1)$, we want to find a point $x$ with $a \leq x \leq b$. Now we compute a sequence $(J_i)_{i \in \mathbb{N}}$ of rational intervals that are open in $[0,1]$ by

$$J_{2i} := [0, q_i), J_{2i+1} := (r_i, 1]$$

and we obtain the interval $I = [a,b] = [0,1] \setminus \bigcup_{i=0}^{\infty} J_i$. With the help of a realizer of $\mathsf{C_I}$, we can find a point $x \in I$, i.e. such that $a \leq x \leq b$. This shows $\mathsf{B_I} \equiv_{\mathrm{sW}} \mathsf{B_I}' \leq_{\mathrm{sW}} \mathsf{C_I}$. The same proof shows $\mathsf{B_I}^- \equiv_{\mathrm{sW}} \mathsf{C_I}^-$.

The reduction $\mathsf{B_I}^+ \leq_{\mathrm{sW}} \mathsf{C_A}$ can be proved analogously $\mathsf{B_I}' \leq_{\mathrm{sW}} \mathsf{C_I}$ using the intervals

$$J_{2i} := (-\infty, q_i), J_{2i+1} := (r_i, \infty),$$

where $(\infty, \infty) = \emptyset$. ⊣

We recall that it is known that $\mathsf{B}$ is equivalent to $\mathsf{C} = \widehat{\mathsf{LPO}}$. The map $\mathsf{C}$ can be considered as countable closed choice and it is equivalent to the map $\mathcal{A}_-(\mathbb{N}) \to \{0,1\}^{\mathbb{N}}, A \mapsto \mathrm{cf}_A$.

PROPOSITION 3.7 (Countable closed choice). $\mathsf{B} \equiv_{\mathrm{W}} \mathsf{C} \equiv_{\mathrm{W}} \widehat{\mathsf{LPO}}$.

PROOF. It is well known that translating a $\rho_<$–name of a real number into a $\rho$–name is equivalent to $\mathsf{C}$, this has originally been proved in [50] (see also Exercise 8.2.12 in [52]). Since $\mathsf{C} = \widehat{\mathsf{LPO}}$ this means $\mathsf{B} \equiv_{\mathrm{W}} \widehat{\mathsf{LPO}}$. ⊣

Now we will show that proper interval choice is reducible to discrete choice. This might be surprising on the first sight and indeed one does not get strict reducibility but one has to exploit direct access to the input in an essential way for this reduction, similarly to the proof of the reduction $\mathsf{C_F} \leq_{\mathrm{W}} \mathsf{C_{\mathbb{N}}}$ in Proposition 3.3.

PROPOSITION 3.8 (Proper interval and discrete choice). $\mathsf{C_I}^- \leq_{\mathrm{W}} \mathsf{C_{\mathbb{N}}}$.

PROOF. We have $\mathsf{C_{\mathbb{N}}} \equiv_{\mathrm{W}} \mathsf{B_F}$ by Proposition 3.3 and $\mathsf{B_I}^- \equiv_{\mathrm{sW}} \mathsf{C_I}^-$ by Proposition 3.6. Thus it suffices to prove $\mathsf{B_I}^- \leq_{\mathrm{W}} \mathsf{B_F}$. Given a name $p$ that encodes an increasing sequence $(q_n)_{n \in \mathbb{N}}$ of rational numbers and a decreasing sequence of rational numbers $(r_n)_{n \in \mathbb{N}}$ with $x := \sup_{n \in \mathbb{N}} q_n < \inf_{n \in \mathbb{N}} r_n =: y$ we want to find a point $z \in \mathbb{R}$ with $x \leq z \leq y$. We inspect the sequences of rational numbers for indexes $n = 0, 1, 2....$ At stage 0 we "'guess" that the center $c_0 := q_0 + \frac{r_0 - q_0}{2}$ is a good choice for a point $z$. If at stage $n$ we have some current choice $c$, then we check whether this choice is supported in this stage. If $q_n > c$ or $r_n < c$, then we dismiss the current guess $c$ and replace it by the new guess $c_n := q_n + \frac{r_n - q_n}{2}$, the center of $[q_n, r_n]$. We record those stages $n$ where we have to dismiss the guess in a sequence $r$, starting with $r(0) = 0$. In those stages where we do not have to dismiss the guess, we just repeat the number of the previous stage that was dismissive in $r$. Altogether, this describes a computation of a function $K$ with $K(p) = r$. Since $x < y$, it is clear that the sequence $r$ of rational numbers is bounded. Using a realizer of $\mathsf{B_F}$ we obtain a corresponding upper bound $m \geq \max \mathrm{range}(r)$. Using this bound $m$ and the original input $p$, we can repeat the computation of guesses of centers starting with $c_0$ until stage $m$. We then use the current guess $c$ at stage $m$ as final choice. The upper bound $m$ guarantees that this choice $z = c$ will never be dismissed in the future and hence satisfies $x \leq z \leq y$. Altogether, this describes how $\mathsf{B_I}^-$ can be reduced to $\mathsf{B_F}$.          ⊣

We note that in the end of the proof it is necessary to repeat the whole computation until stage $m$ and we cannot simply choose the center that is reached at stage $m$. This is because the center of stage $m$ is not necessarily better than the previous centers at some stage $k \leq m$. Figure 3 illustrates the situation (the center $c$ of stage 0 survives stage $m + 1$ but the center of stage $m$ would not survive).

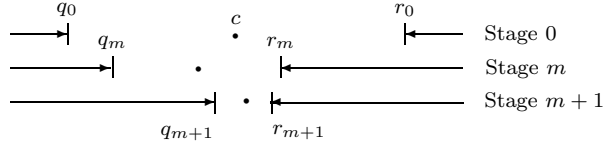

FIGURE 3. Centers of stages $0, m, m + 1$

Now we prove that $\mathsf{LLPO}$ can be reduced to proper interval choice $\mathsf{C_I}^-$ and $\mathsf{LPO}$ can be reduced to discrete choice $\mathsf{C_N}$.

PROPOSITION 3.9 (Omniscience and choice). $\mathsf{LLPO} \leq_{\mathrm{W}} \mathsf{C_I}^-$, $\mathsf{LPO} \leq_{\mathrm{W}} \mathsf{C_N}$.

PROOF. By Propositions 3.3 and 3.6 it suffices to prove the two reductions $\mathsf{LPO} \leq_{\mathrm{W}} \mathsf{B_F}$ and $\mathsf{LLPO} \leq_{\mathrm{W}} \mathsf{B_I}^-$. We prove that $\mathsf{LPO} \leq_{\mathrm{W}} \mathsf{B_F}$ holds. For this purpose we consider a computable function $f : \mathbb{N}^{\mathbb{N}} \to \mathbb{R}_<, p \mapsto x_p$, where

$$x_p := \begin{cases} \min\{n \in \mathbb{N} : p(n) \neq 0\} & \text{if } p \neq 0^{\mathbb{N}} \\ 0 & \text{otherwise} \end{cases}$$

Then we obtain for all $p \neq 0^{\mathbb{N}}$

$$y \in \mathsf{B_F}(x_p) \iff x_p \leq y \iff (\exists n \leq y)\, p(n) \neq 0.$$

Thus, any realizer of $\mathsf{B_F}$ can be used to compute a realizer of $\mathsf{LPO}$, since the realizer of $\mathsf{B_F}$ yields a "search modulus" $y$ that one can use to determine a suitable value for $\mathsf{LPO}$ with direct access to $p$. More precisely, for every $y \in \mathsf{B_F}(x_p)$ we obtain

$$p \neq 0^{\mathbb{N}} \iff p[y + 1] = p(0)...p(y) \neq 0^{y+1}.$$

This means $\mathsf{LPO} \leq_{\mathrm{W}} \mathsf{B_F}$.

We prove $\mathsf{LLPO} \leq_{\mathrm{W}} \mathsf{B_I}^-$. Given a sequence $p \in \mathbb{N}^{\mathbb{N}}$ that is different from 0 in at most one position we search for such a position. As long as $p(n) = 0$ for $n = 0, 1, 2, ..., $, we generate two sequences of values $q_n = -2^{-n}$ and $r_n = 1 + 2^{-n}$. As soon as we find a value $n$ with $p(n) \neq 0$, we continue differently. If $n$ is even, then we continue with the values $r_k = 1/4 + 2^{-k}$ for all $k \geq n$, whereas the $q_k$ are left as described above. If $n$ is odd, then we continue with the values $q_k = 3/4 - 2^{-k}$ for $k \geq n$, whereas the values for $r_k$ are left as described originally. Let $x := \sup_{n \in \mathbb{N}} q_n$ and $y := \inf_{n \in \mathbb{N}} r_n$. Now we obtain

$$[x,y] = \left\{ \begin{array}{ll} [0,1] & \text{if } p(n) = 0 \text{ for all } n \\ [0,1/4] & \text{if } p(2n) \neq 0 \text{ for some } n \\ [3/4,1] & \text{if } p(2n+1) \neq 0 \text{ for some } n \end{array} \right. .$$

Given any realizer of $\mathsf{B_I}^-$, we can now determine a value $z \in \mathsf{B_I}^-(x,y) = [x,y]$. Then we check whether $z \in (1/4, 1]$ or $z \in [0, 3/4)$ and depending which test gives a positive result first, we produce the output $m = 0$ or $m = 1$. In any case $m \in \mathsf{LLPO}(p)$. Thus $\mathsf{LLPO} \leq_{\mathrm{sW}} \mathsf{B_I}^-$. The construction of intervals is illustrated in Figure 4. ⊣

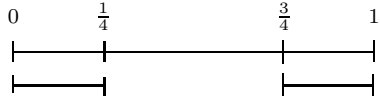

FIGURE 4. Interval construction for $\mathsf{LLPO}$

If we compile all the observations of this section then we obtain the following main result. In fact, we have identified two chains of choice principles that are related in the given way. This clarifies all the positive relations between the corresponding operations given in Figure 1.

THEOREM 3.10 (Choice chains). *We obtain*

1. $\mathsf{LLPO} \leq_{\mathrm{W}} \mathsf{C_I}^- \leq_{\mathrm{W}} \mathsf{C_I} \leq_{\mathrm{W}} \mathsf{C_K} \equiv_{\mathrm{W}} \widehat{\mathsf{LLPO}} \leq_{\mathrm{W}} \mathsf{C_A}$.
2. $\mathsf{LPO} \leq_{\mathrm{W}} \mathsf{C_{\mathbb{N}}} \leq_{\mathrm{W}} \mathsf{B_I}^+ \leq_{\mathrm{W}} \mathsf{C_A} \leq_{\mathrm{W}} \mathsf{C} \equiv_{\mathrm{W}} \widehat{\mathsf{LPO}}$.
3. $\mathsf{LLPO} \leq_{\mathrm{W}} \mathsf{LPO}$, $\mathsf{C_I}^- \leq_{\mathrm{W}} \mathsf{C_{\mathbb{N}}}$, $\mathsf{C_I} \leq_{\mathrm{W}} \mathsf{B_I}^+$.

PROOF. We recall that $\mathsf{C} = \widehat{\mathsf{LPO}}$, $\mathsf{C_K} \equiv_{\mathrm{W}} \widehat{\mathsf{LLPO}}$ by Theorem 2.10 and by Theorem 2.7 we obtain $\mathsf{LLPO} \leq_{\mathrm{W}} \mathsf{LPO}$. The reduction $\mathsf{C_A} \leq_{\mathrm{W}} \mathsf{C}$ since $\mathsf{C_A}$ is $\boldsymbol{\Sigma}_2^0$–computable by Proposition 4.5 of [10] and $\mathsf{C}$ is known to be $\boldsymbol{\Sigma}_2^0$–complete [5]. By Proposition 3.8 we have $\mathsf{C_I}^- \leq_{\mathrm{W}} \mathsf{C_{\mathbb{N}}}$ and by Proposition 3.6 we have $\mathsf{B_I}^+ \leq_{\mathrm{W}} \mathsf{C_A}$. By Proposition 3.9 we have $\mathsf{LLPO} \leq_{\mathrm{W}} \mathsf{C_I}^-$ and $\mathsf{LPO} \leq_{\mathrm{W}} \mathsf{C_{\mathbb{N}}}$.

It is clear that $\mathsf{C_I}^- \leq_{\mathrm{W}} \mathsf{C_I} \leq_{\mathrm{W}} \mathsf{C_K} \leq_{\mathrm{W}} \mathsf{C_A}$ since each operation in this chain is a restriction of the next one. Similarly, $\mathsf{B_I} \leq_{\mathrm{W}} \mathsf{B_I}^+$ since $\mathsf{B_I}$ is a restriction of $\mathsf{B_I}^+$ and $\mathsf{B_F} \leq_{\mathrm{W}} \mathsf{B_I}^+$ since $\mathsf{B_F}(x) = \mathsf{B_I}^+(x, \infty)$, which implies $\mathsf{C_I} \leq_{\mathrm{W}} \mathsf{B_I}^+$ and $\mathsf{C_{\mathbb{N}}} \leq_{\mathrm{W}} \mathsf{B_I}^+$ by Propositions 3.6 and 3.3. ⊣

In the next section we will see that all the given reductions are strict (besides those cases where we have stated equivalence).

Theorem 2.7 and the fact that parallelization is a closure operator by Proposition 2.5 allow us to conclude that the parallelized version of our choice principles are all either equivalent to $\widehat{\mathsf{LLPO}}$ or to $\widehat{\mathsf{LPO}}$ (which are in the relation $\widehat{\mathsf{LLPO}} <_{\mathrm{W}} \widehat{\mathsf{LPO}}$ to each other).

COROLLARY 3.11 (Countable choice principles). *We obtain the following two equivalence classes:* $\widehat{\mathsf{LLPO}} \equiv_{\mathrm{W}} \widehat{\mathsf{C_I^-}} \equiv_{\mathrm{W}} \widehat{\mathsf{C_I}} \equiv_{\mathrm{W}} \widehat{\mathsf{C_K}} <_{\mathrm{W}} \widehat{\mathsf{LPO}} \equiv_{\mathrm{W}} \widehat{\mathsf{C_\mathbb{N}}} \equiv_{\mathrm{W}} \widehat{\mathsf{C_A}}.$

§4. **Separation techniques.** In this section we discuss a number of separation techniques for Weihrauch degrees. These techniques include the Parallelization Principle, the Mind Change Principle, the Computable Invariance Principle, the Low Invariance Principle and the Baire Category Principle. The Mind Change Principle, the Computable and the Low Invariance Principle are all instances of a more general Invariance Principle. The idea is basically to identify a property of realizers that is preserved downwards by Weihrauch reducibility and that distinguishes the two degrees that are to be separated.

The first technique that we mention is the parallelization principle. It is a straightforward consequence of Theorem 2.7 on the omniscience principles and the fact that parallelization is a closure operator.

LEMMA 4.1 (Parallelization principle). *If* $f :\subseteq X \rightrightarrows Y$ *is such that* $\mathsf{LPO} \leq_{\mathrm{W}} f$ *and* $g :\subseteq U \rightrightarrows V$ *is such that* $g \leq_{\mathrm{W}} \widehat{\mathsf{LLPO}}$, *then* $f \not\leq_{\mathrm{W}} g$.

PROOF. Let $\mathsf{LPO} \leq_{\mathrm{W}} f$ and $g \leq_{\mathrm{W}} \widehat{\mathsf{LLPO}}$. Let us assume that $f \leq_{\mathrm{W}} g$ holds. Since parallelization is a closure operator by Proposition 2.5 we obtain

$$\widehat{\mathsf{LPO}} \leq_{\mathrm{W}} \widehat{f} \leq_{\mathrm{W}} \widehat{g} \leq_{\mathrm{W}} \widehat{\mathsf{LLPO}}$$

in contradiction to Theorem 2.7.                                                    ⊣

As an example of an application of the parallelization principle we mention that it implies that co-finite choice $\mathsf{C_F}$, discrete choice $\mathsf{C_\mathbb{N}}$ and closed choice $\mathsf{C_A}$ are all not reducible to compact choice $\mathsf{C_K}$.

COROLLARY 4.2. $\mathsf{C_\mathbb{N}} \equiv_{\mathrm{W}} \mathsf{C_F} \not\leq_{\mathrm{W}} \mathsf{C_K}$ *and* $\mathsf{C_A} \not\leq_{\mathrm{W}} \mathsf{C_K}$.

Several other negative results such as $\mathsf{B_I^+} \not\leq_{\mathrm{W}} \mathsf{C_I}$ follow in the same way, but we do not give a complete list. Now we want to consider the Invariance Principles and as a preparation we define computations with mind changes[1].

DEFINITION 4.3 (Mind changes). Let $n \in \mathbb{N}$. We say that a function $F :\subseteq \mathbb{N}^\mathbb{N} \to \mathbb{N}^\mathbb{N}$ is *computable with at most* $n$ *mind changes*, if there is a limit Turing machine such that for a given input $p$ the machine produces output $F(p)$ in the long run with a two-way output, but such that the machine moves the head on the output tape backwards at most $n$ many times (but each time for an arbitrary finite number of positions) and otherwise the machine operates under one-way output conditions.

---

[1] The concept of mind changes seems to be related to the concept of a topological level of a function, as studied by Hertling [23] and Pauly [39].

In fact, those functions that can be computed with 0 mind changes are exactly the computable functions. Now we can combine all our Invariance Principles in a single result. We recall that a point $p \in \mathbb{N}^{\mathbb{N}}$ is called *low*, if $p' \leq_{\mathrm{T}} \emptyset'$ (we can identify $p$ with graph($p$) for all purposes regarding jumps and Turing reducibility).

LEMMA 4.4 (Invariance Principles). *Let $X, Y, U$ and $V$ be represented spaces and let $f :\subseteq X \rightrightarrows Y$ and $g :\subseteq U \rightrightarrows V$ be multi-valued functions such that $f \leq_{\mathrm{W}} g$. Let $n \in \mathbb{N}$.*

1. *(Mind Change Principle) If $g$ has a realizer that can be computed with at most $n$ mind changes, then $f$ has a realizer that can be computed with at most $n$ mind changes.*
2. *(Computable Invariance Principle) If $g$ has a realizer that maps computable inputs to computable outputs, then $f$ has a realizer that maps computable inputs to computable outputs.*
3. *(Low Invariance Principle) If $g$ has a realizer that maps computable inputs to low outputs, then $f$ has a realizer that maps computable inputs to low outputs.*

PROOF. If $f \leq_{\mathrm{W}} g$, then there are computable functions $H, K :\subseteq \mathbb{N}^{\mathbb{N}} \to \mathbb{N}^{\mathbb{N}}$ such that $H\langle \mathrm{id}, GK \rangle$ is a realizer of $f$ for any realizer $G$ of $g$.

1. We prove that if $G$ is computable with at most $n$ mind changes, then $H\langle \mathrm{id}, GK \rangle$ can also be computed with at most $n$ mind changes. Firstly, $GK$ can be computed with at most $n$ mind changes just by composing the corresponding machines with each other. The machine for $K$ requires no mind change and the machine for $G$ at most $n$ mind changes, hence the machine for $GK$ requires at most $n$ mind changes. It is straightforward to see that also $\langle \mathrm{id}, GK \rangle$ can be computed with at most $n$ mind changes. Now we still have to see that also $H\langle \mathrm{id}, GK \rangle$ can be computed with at most $n$ mind changes. Once again we can just compose the machines. In case that the first machine that computes $\langle \mathrm{id}, GK \rangle$ performs a mind change and moves the head back on the output tape, we can just move the head of the machine for $H$ back into the initial state and start the computation of $H$ from scratch. This will happen at most $n$ many times and hence the composition can be computed with at most $n$ mind changes. One should note that the simulation happens within one single machine where the output of the first machine is only simulated and can be revised within the simulation as often as required. Therefore, if $g$ has a realizer $G$ that is computable with at most $n$ mind changes, then also $f$ has a realizer $H\langle \mathrm{id}, GK \rangle$ that is computable with at most $n$ mind changes.

2. If $G$ is a realizer of $g$ such that $G(p)$ is computable for any computable $p$, then $H\langle p, GK(p) \rangle$ is computable for any computable $p$ and hence $f$ has a realizer that maps computable inputs to computable outputs.

3. If $G$ is a realizer of $g$ such that $G(p)$ is low for any computable $p$, then $H\langle p, GK(p) \rangle$ is low for any computable $p$, because

$$(H\langle \mathrm{p}, GK(p) \rangle)' \leq_{\mathrm{T}} (GK(p))' \leq_{\mathrm{T}} \emptyset'$$

for any computable $p$. ⊣

We start with applying the Mind Change Principle. It can be used to separate LPO and LLPO from all Weihrauch degrees above them that we have considered.

PROPOSITION 4.5. $\mathsf{LPO} <_W \mathsf{B_F}$ *and* $\mathsf{LLPO} <_W \mathsf{B_I}^-$.

PROOF. The positive statements $\mathsf{LPO} \leq_W \mathsf{B_F}$ and $\mathsf{LLPO} \leq_W \mathsf{B_I}^-$ have been proved in Theorem 3.10. The Mind Change Principle in Lemma 4.4 implies the negative statements $\mathsf{B_F} \not\leq_W \mathsf{LPO}$ and $\mathsf{B_I}^- \not\leq_W \mathsf{LLPO}$.

Obviously, LPO can be computed with one mind change. Given $p \in \mathrm{dom}(\mathsf{LPO})$, a Turing machine can just bet on output 1 and inspect the input until a value $n$ with $p(n) = 0$ is found. In this case the machine has to change its mind and revise the output to 0. No further revisions are required, since this is the final result. Similarly, LLPO has a realizer that can be computed with one mind change.

It suffices to prove now that $\mathsf{B_F}$ and $\mathsf{B_I}^-$ have no realizers that can be computed with at most one mind change. We recall that we can assume that the range of $\mathsf{B_F}$ is $\mathbb{N}$. Let us assume that $M$ is a Turing machine that computes a realizer of $\mathsf{B_F}$ with only one mind change. Given a sequence of rational numbers $(q_i)_{i \in \mathbb{N}}$ with $\sup_{i \in \mathbb{N}} q_i = 0$ the machine will eventually produce an output $n \geq 0$. Until this time step the machine has only seen a prefix $q_0, ..., q_n$ of the input sequence $(q_i)_{i \in \mathbb{N}}$ and this sequence can be extended to a sequence $q_0, ..., q_n, n+1, n+1, ...$ that has supremum $n + 1$ which is strictly greater than $n$. Since the machine $M$ will start to operate on this new input in the same way, it will have to perform a mind change after a while and revise the output to some value $k \geq n + 1$. After a certain number of time steps such an output $k$ is produced and once again one can modify the input to have a larger supremum than $k$ and this forces the machine to make another mind change at a later stage. Thus, the assumption that $M$ makes only one mind change is not correct.

One can prove in a similar way that $\mathsf{B_I}^-$ has no realizer that can be computed with at most one mind change. Let us assume that $M$ is a machine that computes a realizer of $\mathsf{B_I}^-$. Given two sequences $(q_i)_{i \in \mathbb{N}}$ and $(r_i)_{i \in \mathbb{N}}$ of rational numbers such that $\sup_{i \in \mathbb{N}} q_i = a < b = \inf_{i \in \mathbb{N}} r_i$, eventually at a time step $t$ the machine $M$ has to produce a result $x$ with $a \leq x \leq b$ with precision $\varepsilon < (b-a)/2$. Now one can modify one of the input sequences to become $(q'_i)_{i \in \mathbb{N}}$ and $(r'_i)_{i \in \mathbb{N}}$ such that these sequences coincide with $(q_i)_{i \in \mathbb{N}}$ and $(r_i)_{i \in \mathbb{N}}$ on a prefix of length at least $t$, respectively, but such that $(x - \varepsilon, x + \varepsilon) \cap [a', b'] = \emptyset$, where $\sup_{i \in \mathbb{N}} q'_i = a' < b' = \inf_{i \in \mathbb{N}} r'_i$. This forces machine $M$ to make a mind change and the process can be repeated as above. Thus, $\mathsf{B_I}^-$ has no realizer that can be computed with only one mind change. ⊣

Now we want to illustrate an application of the Computable Invariance Principle. In particular, it implies that compact choice $\mathsf{C_K}$ cannot be reduced to interval choice $\mathsf{C_I}$ and closed choice $\mathsf{C_A}$ cannot be reduced to discrete choice $\mathsf{C_N}$.

PROPOSITION 4.6. $\mathsf{C_I} \equiv_W \mathsf{B_I} <_W \mathsf{C_K}$ *and* $\mathsf{C_N} \leq_W \mathsf{B_I}^+ <_W \mathsf{C_A}$.

PROOF. The positive statements $\mathsf{C_I} \equiv_W \mathsf{B_I} \leq_W \mathsf{C_K}$ and $\mathsf{C_N} \leq_W \mathsf{B_I}^+ \leq_W \mathsf{C_A}$ have been proved in Theorem 3.10 and Proposition 3.6. The Computable Invariance Principle in Lemma 4.4 implies the negative statements $\mathsf{C_K} \not\leq_W \mathsf{B_I}$ and $\mathsf{C_A} \not\leq_W \mathsf{B_I}^+$. It is straightforward to see that $\mathsf{B_I}$ and $\mathsf{B_I}^+$ have realizers that map computable inputs to computable outputs. This is simply because any non-empty interval $[a, b]$ with a left computable $a$ and a right computable $b$ has a computable member. This can either be a rational number, if $a < b$ or in case $a = b$ it is $a = b$ itself. On the other hand, $\mathsf{C_A}$ and $\mathsf{C_K}$ do not admit realizers that map all computable inputs to computable outputs. This is because there are examples such as Kreisel and Lacombe's set $K \subseteq [0, 1]$ that is a non-empty co-c.e. compact set, which does not contain any computable points (see [34] or Theorem 6.3.8.2 in [52]). ⊣

In order to illustrate an application of the Low Invariance Principle formulated in Lemma 4.4, we prove a version of the Low Basis Theorem of Jockusch and Soare (see [28] or Theorem V.5.32 in [38]) for real number subsets. Indirectly, it is based on the fact that the real numbers are locally compact.

THEOREM 4.7 (Real Low Basis Theorem). *Let $A \subseteq \mathbb{R}$ be a non-empty co-c.e. set. Then there is some point $x \in A$ which has a low $\rho$–name, i.e. there is a $p \in \mathbb{N}^\mathbb{N}$ such that $p' \leq_\mathrm{T} \emptyset'$ and $\rho(p) = x$.*

PROOF. If $A \subseteq \mathbb{R}$ is non-empty, then there is a rational interval $[a, b]$ such that $B := A \cap [a, b] \neq \emptyset$. For this interval $[a, b]$ there is a total representation $\rho_{[a,b]} : \{0, 1\}^\mathbb{N} \to [a, b]$ of $[a, b]$ that is computably equivalent to $\rho$ restricted to $[a, b]$. This can be seen as follows. One can just use the signed-digit representation $\rho_2 : \{0, 1, -1\}^\mathbb{N} \to [-1, 1]$, defined by

$$\rho_2(p) = \sum_{j=0}^{\infty} p(j) 2^{-j-1},$$

which is known to be computably equivalent to $\rho$ restricted to $[-1, 1]$ (see Theorem 7.2.5 in [52]) and then one uses some standard computable surjection $T : \{0, 1\}^\mathbb{N} \to \{0, 1, -1\}^\mathbb{N}$ and

$$\rho_{[a,b]}(p) := \frac{1}{2}(\rho_2 T(p)(b - a) + a + b).$$

In particular there is a computable function $F :\subseteq \{0, 1\}^\mathbb{N} \to \mathbb{N}^\mathbb{N}$ such that $\rho_{[a,b]}(q) = \rho F(q)$ for all $q \in \{0, 1\}^\mathbb{N}$. Then $C := \rho_{[a,b]}^{-1}(B)$ is a co-c.e. closed subset of Cantor space $\{0, 1\}^\mathbb{N}$ since $B = A \cap [a, b]$ is co-c.e. closed. By the classical Low Basis Theorem of Jockusch and Soare any such set has a member $q \in C$ that is low, i.e. $q' \leq_\mathrm{T} \emptyset'$. Now $\rho F(q) = \rho_{[a,b]}(q) = x \in A$ and $p = F(q)$ satisfies $p \leq_\mathrm{T} q$ and hence $p' \leq_\mathrm{T} q' \leq_\mathrm{T} \emptyset'$. Thus $p$ is low too. ⊣

As a consequence of the Low Basis Theorem and the Low Invariance Principle in Lemma 4.4 we obtain that closed choice is $\mathbf{\Sigma_2^0}$–computable, but not $\mathbf{\Sigma_2^0}$–complete. One can also interpret this such that Countable Closed Choice is not reducible to Closed Choice.

PROPOSITION 4.8 (Incompleteness of closed choice). $\mathsf{C_A} <_\mathrm{W} \mathsf{C}$.

PROOF. The reduction $\mathsf{C_A} \leq_\mathrm{W} \mathsf{C}$ has been proved in Theorem 3.10. We need to prove $\mathsf{C} \not\leq_\mathrm{W} \mathsf{C_A}$. To the contrary, let us assume $\mathsf{C} \leq_\mathrm{W} \mathsf{C_A}$. By the Real Low Basis Theorem 4.7 there is a realizer $F$ of $\mathsf{C_A}$ that selects for any computable $\psi_-$–name $q$ of some $A \subseteq \mathbb{R}$ a $\rho$–name $r = F(q)$ which is low. Then, by the Low Invariance Principle in Lemma 4.4 it follows that $\mathsf{C}$ also has a realizer that has a low output for any computable input $p$. Since $\mathsf{C}$ is a function on Baire space, this means that $\mathsf{C}$ itself has the property that $\mathsf{C}(p)$ is low for any computable $p$. But this is a contradiction to Theorem 8.3 in [5], which shows that there is a computable $p$ such that $\mathsf{C}(p) \equiv_{\mathrm{T}} \emptyset'$ and hence not low.                        ⊣

Finally, we illustrate how the Baire Category Theorem can be used as a separation tool. This is particularly interesting, because, as we will show later, $\mathsf{C}_\mathsf{I}$ represents the Intermediate Value Theorem and $\mathsf{C}_\mathbb{N}$ represents the Baire Category Theorem. Thus, we are using the Baire Category Theorem in order to prove that the Baire Category Theorem does not "prove" the Intermediate Value Theorem. We denote by $\widehat{n} \in \mathbb{N}^\mathbb{N}$ the constant sequence with value $n \in \mathbb{N}$ and we use the representation $\delta_\mathbb{N}$ of $\mathbb{N}$, which is defined by $\delta_\mathbb{N}(\widehat{n}) = n$ for all $n \in \mathbb{N}$ (and undefined for all other inputs).

PROPOSITION 4.9. $\mathsf{C}_\mathsf{I} \not\leq_\mathrm{W} \mathsf{C}_\mathbb{N}$.

PROOF. By Propositions 3.3 and 3.6 it suffices to show $\mathsf{B}_\mathsf{I} \not\leq_\mathrm{W} \mathsf{B}_\mathsf{F}$. Without loss of generality we assume that $\rho_<$–names are strictly increasing sequences of rational numbers that converge to the represented point and that $\rho_>$–names are strictly decreasing converging sequences of rational numbers. As usual, we assume that $\mathsf{B}_\mathsf{F}$ has range $\mathbb{N}$. We mention that the set $C := \{\langle p,q \rangle \in \mathbb{N}^\mathbb{N} : \rho_<(p) \leq \rho_>(q)\}$ is co-c.e. closed. We assume that $\delta$ is $[\rho_<, \rho_>]$ restricted to $C$. Let us assume that $\mathsf{B}_\mathsf{I} \leq_\mathrm{W} \mathsf{B}_\mathsf{F}$. Then there are computable functions $H, K$ such that $H\langle \mathrm{id}, GK \rangle$ is a $(\delta, \rho)$–realizer of $\mathsf{B}_\mathsf{I}$ for any $(\rho_<, \delta_\mathbb{N})$–realizer $G$ of $\mathsf{B}_\mathsf{F}$. We can assume that $\mathrm{dom}(K) = C$ (otherwise we restrict $K$ to $C$). Now we note that

$$P_n := \{p \in \mathbb{N}^\mathbb{N} : \rho_< K(p) \leq n\}$$

is closed in $\mathrm{dom}(K) = C$ for all $n \in \mathbb{N}$ since $K$ is continuous. Moreover, $\bigcup_{n=0}^\infty P_n = C$ and by the Baire Category Theorem there must be some $n \in \mathbb{N}$ and some $w \in \mathbb{N}^*$ such that $\emptyset \neq w\mathbb{N}^\mathbb{N} \cap C \subseteq P_n$ since $C$ is a complete metric space and all the $P_n$ are closed in $C$. Let us now fix some realizer $G$ of $\mathsf{B}_\mathsf{F}$ with $\delta_\mathbb{N} G(p) = \max\{n, \lceil \rho_<(p) \rceil\}$. Here $\lceil x \rceil := \min\{z \in \mathbb{Z} : z \geq x\}$ for any real number $x \in \mathbb{R}$. It is clear that for this realizer

$$P_n = \{p \in \mathbb{N}^\mathbb{N} : \rho_\mathbb{N} K(p) \leq n\} = \{p \in \mathbb{N}^\mathbb{N} : \delta_\mathbb{N} GK(p) = n\}.$$

Without loss of generality, we can assume that $w$ is a prefix of a $\delta$–name that is long enough to determine some interval $[a, b]$. We recall that a name $p$ with respect to $\delta = [\rho_<, \rho_>]|_C$ encodes two sequences $(q_n)_{n \in \mathbb{N}}$ and $(r_n)_{n \in \mathbb{N}}$ of rational numbers, which are strictly increasing and decreasing, respectively, such that $\sup_{n \in \mathbb{N}} q_n \leq \inf_{n \in \mathbb{N}} r_n$. Hence, we assume that $w$ is long enough to determine a prefix $q_0, ..., q_i$ with $a = \max\{q_0, ..., q_i\}$ and $r_0, ..., r_i$ with $b = \min\{r_0, ..., r_i\}$. In particular, $a < b$ and there is some $\delta$–name $p$ with $w \sqsubseteq p$ such that $I := \mathsf{B}_\mathsf{I}\delta(p) = [a', b']$ with $a < a' < b' < b$. Then $x := \rho H\langle p, \widehat{n} \rangle \in I$. Therefore $x \in [a', b')$ or $x \in (a', b']$. Suppose for example that $x \in (a', b']$. Take any open interval $J$ such that $x \in J$ and $\inf(J) > a'$. By continuity of $\rho H$ there is some word $v$ with $w \sqsubseteq v \sqsubseteq p$ such that $\rho H \langle v\mathbb{N}^\mathbb{N}, \widehat{n} \rangle \subseteq J$. Without loss of generality, $v$ determines some interval $[a'', b'']$ (in the same manner as $w$ determines $[a, b]$), with $a < a'' < a' < b' < b'' < b$. It is easy to see that there is some $p'$ with $w \sqsubseteq v \sqsubseteq p'$ such that $\mathsf{B}_\mathsf{I}\delta(p') \subseteq (a'', a']$. Therefore $\mathsf{B}_\mathsf{I}\delta(p') \cap J = \emptyset$, and so $\rho H \langle p', GK(p') \rangle \notin J$. But $p' \in w\mathbb{N}^{\mathbb{N}}$ and so $\delta_{\mathbb{N}} GK(p') = n$. Furthermore, $v \sqsubseteq p'$ and so $\rho H \langle p', \widehat{n} \rangle \in J$, which is a contradiction. The case $x \in [a', b']$ can be treated analogously. Figure 5 illustrates the situation.        $\dashv$

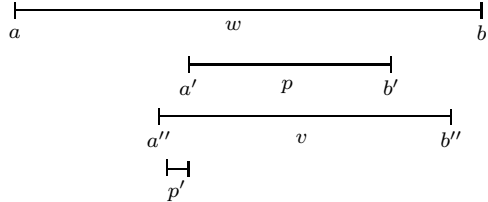

FIGURE 5. Intervals generated by $w, p, v$ and $p'$

Using the previous proposition together with the Parallelization Principle 4.1 we obtain that discrete choice $\mathsf{C}_{\mathbb{N}}$ and interval choice $\mathsf{C}_{\mathsf{I}}$ are incomparable.

COROLLARY 4.10. $\mathsf{C}_{\mathbb{N}} \mid_{\mathrm{W}} \mathsf{C}_{\mathsf{I}}$.

This has a number of negative consequences merely by transitivity, summarized in the following result. In particular, interval choice is not reducible to proper interval choice.

COROLLARY 4.11. $\mathsf{C}_{\mathsf{I}} \not\leq_{\mathrm{W}} \mathsf{C}_{\mathsf{I}}^-$ *and* $\mathsf{B}_{\mathsf{I}}^+ \not\leq_{\mathrm{W}} \mathsf{C}_{\mathbb{N}}$.

The contrary assumptions would lead to contradictions to the results above. Altogether, we can now strengthen the statement of Theorem 3.10 because all the reductions (besides the equivalences, of course) are now known to be strict.

COROLLARY 4.12 (Choice hierarchies). *We obtain*

1. $\mathsf{LLPO} <_{\mathrm{W}} \mathsf{C}_{\mathsf{I}}^- <_{\mathrm{W}} \mathsf{C}_{\mathsf{I}} <_{\mathrm{W}} \mathsf{C}_{\mathsf{K}} \equiv_{\mathrm{W}} \widehat{\mathsf{LLPO}} \leq_{\mathrm{W}} \mathsf{C}_{\mathsf{A}}$.
2. $\mathsf{LPO} <_{\mathrm{W}} \mathsf{C}_{\mathbb{N}} <_{\mathrm{W}} \mathsf{B}_{\mathsf{I}}^+ <_{\mathrm{W}} \mathsf{C}_{\mathsf{A}} <_{\mathrm{W}} \mathsf{C} \equiv_{\mathrm{W}} \widehat{\mathsf{LPO}}$.
3. $\mathsf{LLPO} <_{\mathrm{W}} \mathsf{LPO}$, $\mathsf{C}_{\mathsf{I}}^- <_{\mathrm{W}} \mathsf{C}_{\mathbb{N}}$, $\mathsf{C}_{\mathsf{I}} <_{\mathrm{W}} \mathsf{B}_{\mathsf{I}}^+$.

In particular, all the arrows between the Weihrauch degrees in the diagram in Figure 1 that are mentioned here are now proved to be correct and it is clear that no additional arrows can be added (besides those that follow by transitivity).

The reader can verify that negative results can be relativized, if required. For instance, the Parallelization Principle 4.1 is true in a purely topological version. Similarly, the Real Low Basis Theorem 4.7 and other results can be relativized.

§5. **Discrete Choice and the Baire Category Theorem.** In this section we want to classify the Weihrauch degree of the Baire Category Theorem and some core theorems from functional analysis such as the Banach Inverse Mapping Theorem, the Open Mapping Theorem, the Closed Graph Theorem and the Uniform Boundedness Theorem. It is clear that these theorems are closely related and in fact all these theorems are typically directly derived from the Baire Category Theorem or from each other. The computable content of these theorems has been studied in [3, 4, 6, 9] and we will essentially use results obtained in these sources here.

We start to discuss the Baire Category Theorem. As explained in the introduction, it depends on the exact version of the Baire Category Theorem which Weihrauch degree the theorem has. If one formalizes the version $\mathsf{BCT}_0$, as mentioned in the introduction, then one obtains a fully computable theorem (this has been proved in [3]) and we will not discuss this version any further here. The contrapositive version $\mathsf{BCT}$, however, is interesting for us.

**Definition 5.1** (Baire Category Theorem). *Let $X$ be a complete computable metric space. We consider the operation*

$$\mathsf{BCT}_X :\subseteq \mathcal{A}_-(X)^{\mathbb{N}} \rightrightarrows \mathbb{N}, (A_i)_{i \in \mathbb{N}} \mapsto \{n \in \mathbb{N} : A_n \text{ has non-empty interior}\},$$

*where* $\operatorname{dom}(\mathsf{BCT}_X) = \{(A_i)_{i \in \mathbb{N}} : \bigcup_{i=0}^{\infty} A_i = X\}$.

We note that $\mathsf{BCT}_X$ is well-defined by the Baire Category Theorem. We mention a few facts that are known about $\mathsf{BCT}$ (see [3]):

1. $\mathsf{BCT}_X$ is discontinuous and hence non-computable.
2. $\mathsf{BCT}_X$ is obviously non-uniformly computable in the sense that computable inputs are mapped to computable outputs (this is because it has natural number output).
3. $\mathsf{BCT}_X$ admits a non-computable sequential counterexample.

There are other possible ways to formalize this theorem. For instance, one could consider the space as a part of the input, which would however require the discussion of representations of spaces. This approach has been adopted for the Hahn-Banach Theorem in [20]. We restrict the discussion here to a fixed space $X$ and in this case the result is even stronger, since the Baire Category Theorem is equivalent to discrete choice for any fixed non-empty complete computable metric space $X$.

**Theorem 5.2** (Baire Category Theorem). *Let $X$ be a non-empty complete computable metric space. Then* $\mathsf{BCT}_X \equiv_{\mathrm{W}} \mathsf{C}_{\mathbb{N}}$.

Proof. We prove $\mathsf{BCT}_X \leq_{\mathrm{sW}} \mathsf{C}_{\mathbb{N}}$. Given a sequence $(A_i)_{i \in \mathbb{N}}$ of closed sets with $X = \bigcup_{i=0}^{\infty} A_i$ we want to find some $n$ such that $A_n$ has non-empty interior. Let $(B_i)_{i \in \mathbb{N}}$ be an effective enumeration of the open rational balls of $X$. The set

$$P = \{\langle n, m \rangle \in \mathbb{N} : (X \setminus A_n) \cap B_m \neq \emptyset \text{ or } B_m = \emptyset\}$$

is c.e. in and hence we can compute a $p \in \mathbb{N}^{\mathbb{N}}$ such that $\operatorname{range}(p) - 1 = P$. Then $P = \mathbb{N} \setminus Q$ with

$$Q = \{\langle n, m \rangle \in \mathbb{N} : \emptyset \neq B_m \subseteq A_n\}$$

and by the Baire Category Theorem this set is non-empty. Hence, using a realizer of $\mathsf{C}_{\mathbb{N}}$, we can determine a point $\langle n, m \rangle \in Q = \mathbb{N} \setminus P$ and the component $n$ is the desired result. This proves $\mathsf{BCT}_X \leq_{\mathrm{sW}} \mathsf{C}_{\mathbb{N}}$.

In order to prove $\mathsf{C}_{\mathbb{N}} \leq_{\mathrm{W}} \mathsf{BCT}_X$, it suffices to prove $\mathsf{B}_{\mathsf{F}} \leq_{\mathrm{W}} \mathsf{BCT}_X$ by Proposition 3.3. Given a sequence of rational numbers $(q_i)_{i \in \mathbb{N}}$ which is bounded from above, we want to compute some upper bound $n \in \mathbb{N}$ with $x := \sup_{i \in \mathbb{N}} q_i \leq n$. Given $(q_i)_{i \in \mathbb{N}}$ we can compute a sequence $(A_n)_{n \in \mathbb{N}}$ of closed sets with

$$A_n := \begin{cases} \emptyset & \text{if } (\exists i) \, n < q_i \\ X & \text{otherwise} \end{cases}$$

for all $n \in \mathbb{N}$. Then $\bigcup_{n=0}^{\infty} A_n = X$ and any $n$ such that $A_n$ has non-empty interior satisfies $x \leq n$. Such an $n$ can be determined with a realizer of $\mathsf{BCT}_X$ and hence $\mathsf{C}_\mathbb{N} \leq_\mathrm{W} \mathsf{B}_\mathsf{F} \leq_\mathrm{W} \mathsf{BCT}_X$. ⊣

We mention that a seemingly different versions of the Baire Category Theorem are in fact equivalent to $\mathsf{BCT}_X$ as the same proof above shows:

1. Together with the index $n$ of a set $A_n$ that has non-empty interior we can even determine an index $m$ of a non-empty ball $B_m \subseteq A_n$. This is what is required in practice when the Baire Category Theorem is applied.
2. The Baire Category Theorem could be restricted to monotone sequences $(A_n)_{n \in \mathbb{N}}$, i.e. sequences such that $A_n \subseteq A_{n+1}$ for all $n \in \mathbb{N}$.

As a next theorem we discuss the Banach Inverse Mapping Theorem. If $X$ and $Y$ are computable metric spaces with Cauchy representations $\delta_X$ and $\delta_Y$, respectively, then we denote by $\mathcal{C}(X, Y)$ the set of (relatively) continuous functions $T : X \to Y$ and we represent $\mathcal{C}(X, Y)$ by the canonical function space representation $[\delta_X \to \delta_Y]$ (see [52]). A computable Banach space is just a normed space together with a sequence whose linear span is dense in the space and such that the induced metric space is a complete computable metric space (see [9] for details). The Banach Inverse Mapping Theorem can even be described by a single-valued function in a very natural way. It is just the inversion operator.

DEFINITION 5.3 (Banach Inverse Mapping Theorem). *Let $X, Y$ be computable Banach spaces. Then we define*

$$\mathsf{IMT}_{X,Y} :\subseteq \mathcal{C}(X, Y) \to \mathcal{C}(Y, X), T \mapsto T^{-1}$$

*with* $\mathrm{dom}(\mathsf{IMT}_{X,Y}) := \{T \in \mathcal{C}(X, Y) : T \text{ linear and bijective}\}$.

We note that $\mathsf{IMT}_{X,Y}$ is well-defined by the Banach Inverse Mapping Theorem. We briefly note what is known about $\mathsf{IMT}_{X,Y}$ (see [9]):

1. $\mathsf{IMT}_{X,Y}$ is discontinuous and hence non-computable in general.
2. $\mathsf{IMT}_{X,Y}$ is non-uniformly computable, i.e. the inverse $T^{-1}$ of any computable linear and bijective operator $T : X \to Y$ is computable.
3. $\mathsf{IMT}_{X,Y}$ admits a non-computable sequential counterexample for the Hilbert space $X = Y = \ell_2$.
4. $\mathsf{IMT}_{X,Y}$ is computable for finite-dimensional $X, Y$.

This means, in particular, that unlike the situation with the Baire Category Theorem it is the case that the Weihrauch degree of the Banach Inverse Mapping Theorem depends on the underlying spaces. The worst case, however, is already achieved for the computable infinite-dimensional Hilbert space $X = Y = \ell_2$.

THEOREM 5.4 (Banach Inverse Mapping Theorem). *Let $X, Y$ be computable Banach spaces. Then* $\mathsf{IMT}_{X,Y} \leq_\mathrm{W} \mathsf{C}_\mathbb{N} \equiv_\mathrm{W} \mathsf{IMT}_{\ell_2, \ell_2}$.

PROOF. By Proposition 3.3 it is sufficient to prove the claim for $\mathsf{B}_\mathsf{F}$ instead of $\mathsf{C}_\mathbb{N}$. We prove $\mathsf{IMT}_{X,Y} \leq_\mathrm{W} \mathsf{B}_\mathsf{F}$. In Theorem 4.9 of [9] it has been proved that the map

$$\iota :\subseteq \mathcal{C}(X, Y) \times \mathbb{R} \to \mathcal{C}(Y, X), (T, s) \mapsto T^{-1}$$

with $\mathrm{dom}(\iota) = \{(T, s) : T \text{ linear and bijective and } ||T^{-1}|| < s\}$ is computable. On the other hand, it is clear that $||T^{-1}|| = 1/(\inf_{||x||=1} ||Tx||)$ (see Exercise 5.14

in [49]) and a $\rho_<$–name of this number can be computed, for a given $T$. Thus, given a realizer of $\mathsf{B_F}$, we can actually compute an upper bound of $||T^{-1}||$ and together with the input $T$, we can compute $T^{-1}$ using $\iota$.

Now we assume that $X = Y = \ell_2$ and we prove $\mathsf{B_F} \leq_W \mathsf{IMT}_{\ell_2,\ell_2}$. In Proposition 4.7 of [9] it has been proved that there is a computable operation $\tau :\subseteq \mathbb{R}_> \rightrightarrows \mathcal{C}(\ell_2, \ell_2)$ such that for any $a \in \mathbb{R}_>$ with $a \in (0,1]$ there exists some $T_a \in \tau(a)$ and all such $T_a : \ell_2 \to \ell_2$ are linear, bounded, bijective and satisfy $||T_a^{-1}|| = 1/a$. Thus, given some $x \in \mathbb{R}_<$ that without loss of generality satisfies $x \geq 1$, we can compute some $T_{1/x} \in \tau(1/x)$ and using some realizer of $\mathsf{IMT}_{\ell_2,\ell_2}$ we obtain $T_{1/x}^{-1}$. This allows to compute some upper bound $M \geq ||T_{1/x}^{-1}|| = x$ by Theorem 5.1 of [6]. Hence $\mathsf{B_F} \leq_W \mathsf{IMT}_{\ell_2,\ell_2}$.                              $\dashv$

Roughly speaking, the Banach Inverse Mapping Theorem is almost computable. The only extra information that cannot be extracted from the given operator $T$ is the norm $||T^{-1}||$ of its inverse. In fact, the norm can be computed from below and one upper bound suffices as extra information to obtain the inverse. It is exactly this upper bound that can be provided by $\mathsf{B_F} \equiv_W \mathsf{C_N}$. Vice versa, on $\ell_2$ the inverse operator can be used to compute such upper bounds. This situation is practically the same for the Open Mapping Theorem.

DEFINITION 5.5 (Open Mapping Theorem). *Let $X, Y$ be computable Banach spaces. Then we define*

$$\mathsf{OMT}_{X,Y} :\subseteq \mathcal{C}(X,Y) \to \mathcal{C}(\mathcal{O}(X), \mathcal{O}(Y)), T \mapsto (U \mapsto T(U))$$

*with* $\mathrm{dom}(\mathsf{OMT}_{X,Y}) := \{T \in \mathcal{C}(X,Y) : T \text{ linear and surjective}\}$.

The operation $\mathsf{OMT}_{X,Y}$ is well-defined by the classical Open Mapping Theorem. Here $\mathcal{O}(Z)$ denotes the set of open subsets $U \subseteq Z$ and it is represented with the positive information representation $\vartheta$. That is, a name of an open set $U$ is a list of rational open balls $B(x_i, r_i)$ whose union is $U$ (see [13] for details). We just mention the result without proof here (the reader can find all required ingredients for the proof in [9]).

THEOREM 5.6 (Open Mapping Theorem). *Let $X, Y$ be computable Banach spaces. Then* $\mathsf{OMT}_{X,Y} \leq_W \mathsf{C_N} \equiv_W \mathsf{OMT}_{\ell_2,\ell_2}$.

Another interesting theorem is the Closed Graph Theorem. The classical theorem states that a linear operator with a closed graph on Banach spaces is bounded. We can effectivize this result using positive information on the graph.

DEFINITION 5.7 (Closed Graph Theorem). *Let $X, Y$ be computable Banach spaces. Then we define*

$$\mathsf{CGT}_{X,Y} :\subseteq \mathcal{A}_+(X \times Y) \to \mathcal{C}(X,Y), \mathrm{graph}(T) \mapsto T$$

*where* $\mathrm{dom}(\mathsf{CGT}_{X,Y})$ *contains all those closed subsets $A \subseteq X \times Y$ for which there is a linear bounded $T : X \to Y$ such that $A = \mathrm{graph}(T)$.*

Here $\mathcal{A}_+(Z)$ denotes the set of closed subsets of $Z$ equipped with the positive information representation $\psi_+$. A name of a closed set $A \subseteq Z$ in this sense is a list of all rational open balls which intersect $A$.

**Theorem 5.8** (Closed Graph Theorem). *Let $X, Y$ be computable Banach spaces. Then* $\mathsf{CGT}_{X,Y} \leq_{\mathrm{W}} \mathsf{C}_{\mathbb{N}} \equiv_{\mathrm{W}} \mathsf{CGT}_{\ell_2, \ell_2}$.

**Proof.** By Proposition 3.3 it is sufficient to prove the claim for $\mathsf{B}_{\mathsf{F}}$ instead of $\mathsf{C}_{\mathbb{N}}$. By Theorem 4.3(2) in [4] the representation $[\delta_X \to \delta_Y]$ of linear bounded operators $T : X \to Y$ is equivalent to the representation that represents $T$ by positive information on $\mathrm{graph}(T)$ and an upper bound on the operator norm $||T||$. Having positive information on $\mathrm{graph}(T)$ allows to compute a dense sequence within $\mathrm{graph}(T)$ and hence lower bounds of $||T||$. With the help of $\mathsf{B}_{\mathsf{F}}$ one can obtain the required upper bound. This proves $\mathsf{CGT}_{X,Y} \leq_{\mathrm{W}} \mathsf{B}_{\mathsf{F}}$.

Now let $X = Y = \ell_2$. We use Example 5.3 of [6]. Given a real number $a \in \mathbb{R}_<$ by a sequence $(a_i)_{i \in \mathbb{N}}$ of rational numbers with $\sup_{i \in \mathbb{N}} a_i = a$, we compute positive information on $\mathrm{graph}(T)$ of the operator $T : \ell_2 \to \ell_2$ with $T(x_i)_{i \in \mathbb{N}} := (a_i x_i)_{i \in \mathbb{N}}$. Using a realizer of $\mathsf{CGT}_{X,Y}$ we obtain the operator $T$ and hence we can compute some upper bound $M \geq ||T|| = a$ by Theorem 5.1 of [6]. Hence $M \in \mathsf{B}_{\mathsf{F}}(a)$, which proves $\mathsf{B}_{\mathsf{F}} \leq_{\mathrm{W}} \mathsf{CGT}_{\ell_2, \ell_2}$. $\dashv$

Analogously to the Banach Inverse Mapping Theorem and the Open Mapping Theorem it is known that the Closed Graph Theorem is computable for finite-dimensional $X, Y$.

As a final example in this section we discuss the Uniform Boundedness Theorem. In case of the Uniform Boundedness Theorem the relation to $\mathsf{B}_{\mathsf{F}}$ is even more visible that for the other results, since the statement of the theorem is that certain upper bounds exist. We formalize the theorem as follows.

**Definition 5.9** (Uniform Boundedness Theorem). *Let $X, Y$ be a computable Banach spaces. Then we define*

$$\mathsf{UBT}_{X,Y} :\subseteq \mathcal{C}(X,Y)^{\mathbb{N}} \rightrightarrows \mathbb{N}, (T_i)_{i \in \mathbb{N}} \mapsto \{M \in \mathbb{N} : \sup_{i \in \mathbb{N}} ||T_i|| \leq M\}$$

*where $\mathrm{dom}(\mathsf{UBT}_{X,Y})$ contains all sequences $(T_i)_{i \in \mathbb{N}}$ of linear and bounded operators $T_i : X \to Y$ such that $\{||T_i x|| : i \in \mathbb{N}\}$ is bounded for each $x \in X$.*

It follows from the classical Uniform Boundedness Theorem that the operation $\mathsf{UBT}_{X,Y}$ is well-defined. Now we can prove the following equivalence. We mention that similarly as for the Baire Category Theorem it turns out that we obtain the worst case for the Uniform Boundedness Theorem for arbitrary computable Banach spaces $X \neq \{0\}$.

**Theorem 5.10** (Uniform Boundedness Theorem). *Let $X, Y$ be computable Banach spaces different from $\{0\}$. Then* $\mathsf{UBT}_{X,Y} \equiv_{\mathrm{W}} \mathsf{C}_{\mathbb{N}}$.

**Proof.** By Proposition 3.3 it is sufficient to prove the claim for $\mathsf{B}_{\mathsf{F}}$ instead of $\mathsf{C}_{\mathbb{N}}$. In Corollary 5.5 of [6] it has been proved that the uniform bound map

$$U :\subseteq \mathcal{C}(X,Y)^{\mathbb{N}} \to \mathbb{R}_<, (T_i)_{i \in \mathbb{N}} \mapsto \sup_{i \in \mathbb{N}} ||T_i||$$

is computable with $\mathrm{dom}(U) = \mathrm{dom}(\mathsf{UBT}_{X,Y})$. This directly implies the reduction $\mathsf{UBT}_{X,Y} \leq_{\mathrm{W}} \mathsf{B}_{\mathsf{F}}$ and hence we obtain $\mathsf{UBT}_{X,Y} \leq_{\mathrm{W}} \mathsf{C}_{\mathbb{N}}$ with Proposition 3.3.

Now we prove $\mathsf{B}_{\mathsf{F}} \leq_{\mathrm{W}} \mathsf{UBT}_{X,Y}$, using Remark 5.7 of [6]. Given a real number $a \in \mathbb{R}_<$ by a sequence $(a_i)_{i \in \mathbb{N}}$ of rational numbers with $\sup_{i \in \mathbb{N}} a_i = a$, we compute a sequence of operators $T_i : X \to Y$ such that $T_i e := a_i e'$ for some computable unit-length vectors $e \in X$ and $e' \in Y$ and $T_i x := 0$ for all $x \in X$ that are linearly independent of $e$. Then $||T_i|| = a_i$ for all $i \in \mathbb{N}$ and $\{||T_i x|| : i \in \mathbb{N}\}$ is bounded for each $x \in X$. Now using a realizer of $\mathsf{UBT}_{X,Y}$ an upper bound $M \in \mathsf{UBT}_{X,Y}((T_i)_{i \in \mathbb{N}})$ of $\sup_{i \in \mathbb{N}} ||T_i|| = a$ can be determined. Hence $M \in \mathsf{B_F}(a)$, which proves $\mathsf{B_F} \leq_{\mathrm{W}} \mathsf{UBT}_{X,Y}$.                                                                    $\dashv$

Although classically the Banach Inverse Mapping Theorem, the Open Mapping Theorem, the Closed Graph Theorem and the Uniform Boundedness Theorem can be proved with the Baire Category Theorem, we have not exploited any logical relation between these theorems here, but just the fact that they have the same computational power. A common feature of all the theorems discussed in this section that are equivalent to $\mathsf{C_{\mathbb{N}}}$ are:

1. They are discontinuous and hence non-computable (since $\mathsf{C_{\mathbb{N}}}$ is so).
2. They admit non-uniform computable solutions (since $\mathsf{C_{\mathbb{N}}}$ has a realizer that maps computable inputs to computable outputs).
3. They have $\Delta_2^0$–complete sequential counterexamples (since $\widehat{\mathsf{C_{\mathbb{N}}}} \equiv_{\mathrm{W}} \mathsf{C}$, any realizer maps some computable sequence to some $\Delta_2^0$–complete sequence in the arithmetical hierarchy).

All the properties mentioned here are degree theoretic properties and any theorem equivalent to $\mathsf{C_{\mathbb{N}}}$ will be of the same category. In [3], for instance, it required an explicit construction using a simple set to show that the Baire Category Theorem $\mathsf{BCT}$ has a non-computable sequential counterexample. Now such results can be easily derived from the characterization presented here. In particular, all properties of $\mathsf{BCT}_X$ and $\mathsf{IMT}_{X,Y}$ mentioned in the beginning follow immediately from the classification of the Weihrauch degree of the corresponding theorems, except the observation 4. that $\mathsf{IMT}_{X,Y}$ is computable for finite-dimensional $X, Y$.

**§6. Interval Choice and the Intermediate Value Theorem.** In this section we want to study the Intermediate Value Theorem. In the following definition we specify the multi-valued operation which captures this theorem.

DEFINITION 6.1 (Intermediate Value Theorem). *We define a multi-valued operation* $\mathsf{IVT} :\subseteq \mathcal{C}[0,1] \rightrightarrows [0,1]$ *by*

$$\mathsf{IVT}(f) := \{x \in [0,1] : f(x) = 0\}$$

*and* $\mathrm{dom}(\mathsf{IVT}) := \{f \in \mathcal{C}[0,1] : f(0) \cdot f(1) < 0\}$.

This theorem has been carefully analyzed in computable analysis (see Section 6.3 in [52]) and the main results are:

1. $\mathsf{IVT}$ is discontinuous and hence non-computable.
2. $\mathsf{IVT}$ is non-uniformly computable in the sense that $\mathsf{IVT}(f)$ contains a computable point for any computable $f \in \mathrm{dom}(\mathsf{IVT})$.
3. $\mathsf{IVT}$ is a computable single-valued function restricted to functions with a unique zero.
4. $\mathsf{IVT}$ restricted to the set of functions with nowhere dense zero set is computable.

The first three properties of the Intermediate Value Theorem are actually immediate consequences of the following classification of its Weihrauch degree. We prove that the Intermediate Value Theorem is equivalent to Interval Choice.

The last mentioned observation 4. requires an analysis of the Intermediate Value Theorem that goes beyond classifying its Weihrauch degree (it can be proved with the so-called trisection method, a computable variant of bisection).

THEOREM 6.2 (Intermediate Value Theorem). $\mathsf{IVT} \equiv_{\mathrm{sW}} \mathsf{C_I}$.

PROOF. By Proposition 3.6 and Lemma 3.5 it is sufficient to prove the reductions $\mathsf{B_I}' \leq_{\mathrm{sW}} \mathsf{IVT} \leq_{\mathrm{sW}} \mathsf{B_I}$. First we show $\mathsf{IVT} \leq_{\mathrm{sW}} \mathsf{B_I}$. Given a function $f \in \mathcal{C}[0,1]$ with $f(0) \cdot f(1) < 0$ we want to find a zero, i.e. a point $x \in [0,1]$ with $f(x) = 0$. We determine two sequences of rational numbers $(q_n)_{n \in \mathbb{N}}$ and $(r_n)_{n \in \mathbb{N}}$ as follows. We let $q_0 := 0$, $r_0 := 1$. In step $n+1$ of the computation we assume that all the values up to $q_n$ and $r_n$ are given and then we perform an exhaustive search over all possible pairs $(q_{n+1}, r_{n+1})$ of rational numbers for at most $n$ time steps in order to identify the pair with the smallest difference $|q_{n+1} - r_{n+1}|$ that we can find and such that

$$q_n < q_{n+1} < r_{n+1} < r_n \text{ and } f(q_{n+1}) \cdot f(r_{n+1}) < 0.$$

If we cannot find such $q_{n+1}, r_{n+1}$ in $n$ time steps, then we let $q_{n+1} = q_n$ and $r_{n+1} = r_n$. It is clear that in this way we obtain $y := \sup_{n \in \mathbb{N}} q_n$ and $z := \inf_{n \in \mathbb{N}} r_n$ with $[y, z] \subseteq f^{-1}\{0\}$ and any realizer of $\mathsf{B_I}$ can determine a zero $x \in \mathsf{B_I}(y, z)$ of $f$.

Now we prove $\mathsf{B_I}' \leq_{\mathrm{sW}} \mathsf{IVT}$. Given a strictly increasing sequence $(q_n)_{n \in \mathbb{N}}$ and a strictly decreasing sequence $(r_n)_{n \in \mathbb{N}}$ of rational numbers with $y := \sup_{n \in \mathbb{N}} q_n \leq \inf_{n \in \mathbb{N}} r_n := z$ we want to find a number $x$ with $y \leq x \leq z$. Without loss of generality we can assume that $q_0 > 0$ and $r_0 < 1$. Now we compute a sequence of rational polygons $(f_n)_{n \in \mathbb{N}}$ as follows. The function $f_n$ is the polygon with the following vertices:

$$(0, -1), (q_0, -2^{-1}), (q_1, -2^{-2}), ..., (q_n, -2^{-n-1}),$$
$$(r_n, 2^{-n-1}), (r_{n-1}, 2^{-n}), ..., (r_0, 2^{-1}), (1, 1).$$

Figure 6 illustrates the situation. Then $(f_n)_{n \in \mathbb{N}}$ is a sequence that converges effectively to a continuous function $f$ with $f(0) \cdot f(1) < 0$ and $f^{-1}\{0\} = [y, z]$. Thus, any realizer of $\mathsf{IVT}$ can determine a value $x \in \mathsf{IVT}(f)$ with $y \leq x \leq z$. ⊣

We list some common features of all theorems that are equivalent to $\mathsf{C_I}$. We note that by Corollary 2.11 any uniquely determined solution is already uniformly computable in the input.

1. They are discontinuous and hence non-computable (since $\mathsf{C_I}$ is so).
2. They admit non-uniform computable solutions (since $\mathsf{C_I}$ has a realizer that maps computable inputs to computable outputs).
3. They are uniformly computable under all classical conditions where the solution is uniquely determined (since $\mathsf{C_I}$ is weakly computable).
4. They have limit computable sequential counterexamples (since $\widehat{\mathsf{C_I}} \equiv_{\mathrm{W}} \mathsf{WKL}$).
5. They have sequential solutions of any basis type (since $\widehat{\mathsf{C_I}} \equiv_{\mathrm{W}} \mathsf{WKL}$).

By a basis type we mean any set $B \subseteq \mathbb{N}^{\mathbb{N}}$ that forms a basis for $\Pi^0_1$ subsets of Cantor space $\{0, 1\}^{\mathbb{N}}$ (see [17]), such as the set of low points. The fact 4. means for the Intermediate Value Theorem that there exists a computable sequence $(f_n)_{n \in \mathbb{N}}$ of continuous functions $f_n : [0,1] \to \mathbb{R}$ with $f_n(0) \cdot f_n(1) < 0$ such that any sequence $(x_n)_{n \in \mathbb{N}}$ of reals with $f_n(x_n) = 0$ is non-computable. A direct proof of this fact by itself requires almost two pages (see Example 8a by Pour-El and Richards in [40]), whereas here it is a simple corollary of the classification of the degree of the Intermediate Value Theorem.

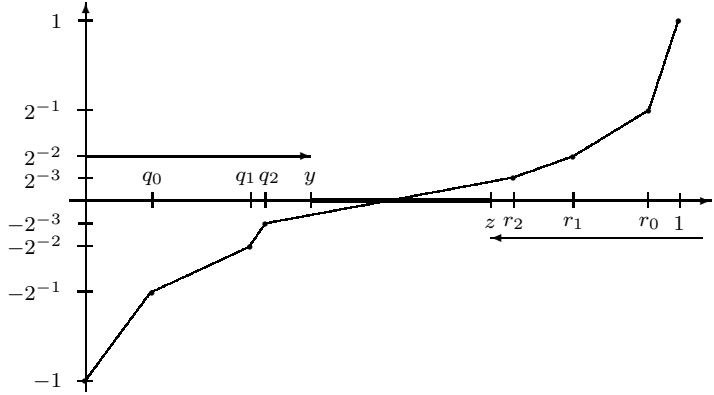

FIGURE 6. Polygon $f_3$

As an observation aside we characterize the image of the set of computable functions under $\mathsf{IVT}$, i.e. the class of zero sets of computable functions that change their sign. We recall that the co-c.e. closed sets $A \subseteq [0,1]$ are exactly those for which there exists a computable function $f : [0,1] \to \mathbb{R}$ with $f^{-1}\{0\} = A$ (see [14, 52]). Somewhat surprisingly, we can attach an arbitrary computable point to an arbitrary co-c.e. closed set to obtain one that appears as a zero set of a function with changing sign.

PROPOSITION 6.3. *Let $A \subseteq (0,1)$ be a non-empty set. Then the following are equivalent:*

1. $A = f^{-1}\{0\}$ *for some computable $f : [0,1] \to \mathbb{R}$ with $f(0) \cdot f(1) < 0$.*
2. $A = B \cup \{x\}$ *for some set $B \subseteq (0,1)$ that is co-c.e. closed in $[0,1]$ and for some computable point $x \in (0,1)$.*
3. $A$ *is co-c.e. closed in $[0,1]$ and contains a computable point $x \in (0,1)$.*

PROOF. We show that 1. implies 3. It is clear that any set $A = f^{-1}\{0\}$ is co-c.e. closed and it follows from the previous Theorem 6.2 that $A$ contains a computable point considering $\mathsf{IVT} \equiv_{\mathrm{W}} \mathsf{B}_2$ (as mentioned above, this is well-known, see [52]). It is clear that 3. implies 2. We show that 2. implies 1. Now let $B \subseteq (0,1)$ be co-c.e. closed in $[0,1]$ and let $x \in (0,1)$ be computable and let $A = B \cup \{x\}$. Then there exists a computable function $g : [0,1] \to \mathbb{R}$ such that $g^{-1}\{0\} = B$. Without loss of generality we can assume $g(0) > 0$ and $g(1) > 0$, since we can replace $g$ by $|g|$ otherwise. Now we define $h(y) := y - x$ for all $y \in [0,1]$ and $f := g \cdot h$ and we get a computable function $f$ with

$$f^{-1}\{0\} = g^{-1}\{0\} \cup h^{-1}\{0\} = B \cup \{x\} = A$$

and $f(0) \cdot f(1) = h(0) \cdot h(1) < 0$.                                    ⊣

**§7. Compact Choice and the Hahn-Banach Theorem.** The Hahn-Banach Theorem has been studied in detail in [8] and [20], see also [12]. We just briefly summarize the known results, without formalizing the corresponding multi-valued function HBT here:

1. HBT is discontinuous and hence non-computable.
2. HBT is computable under all classical conditions that guarantee unique existence of the solution, such as uniform convexity of the dual space, which holds for instance for Hilbert spaces.
3. HBT is not non-uniformly computable in general, not even for finite-dimensional spaces.

The following theorem follows essentially from results in [20], see also [12]. The formalization chosen for this result includes the underlying space as a part of the input data.

THEOREM 7.1 (Hahn-Banach Theorem). $\mathsf{HBT} \equiv_W \mathsf{C_K}$.

Common features of all theorems equivalent to $\mathsf{C_K}$ are:

1. They are discontinuous and hence non-computable (since $\mathsf{C_K}$ is so).
2. They are uniformly computable under all classical conditions where the solution is uniquely determined (since $\mathsf{C_K}$ is weakly computable).
3. They have non-uniform solutions of any basis type (since we have that $\mathsf{C_K} \equiv_W \mathsf{WKL}$).
4. They have limit computable counterexamples (since we have that $\mathsf{C_K} \equiv_W \mathsf{WKL}$).

**§8. Metatheorems and Applications.** In this section we want to discuss a number of metatheorems that allow some conclusions on the status of theorems merely regarding the logical form of these theorems. Essentially, we are trying to identify the computational status of $\Pi_2$–theorems, i.e. theorems of the form

$$(\forall x \in X)(\exists y \in Y)\ (x, y) \in A,$$

where depending on the properties of $Y$ and $A$ automatically certain computable versions of realizers of these theorems exist. In many cases this allows to get some upper bound on the Weihrauch degree of the corresponding theorem straightforwardly. We will discuss these metatheorems together with some characteristic examples. Before we formulate the first metatheorem, we briefly mention open choice.

DEFINITION 8.1 (Open choice). Let $X$ be a computable metric space. The multi-valued operation

$$\mathsf{C}_{\mathcal{O}(X)} :\subseteq \mathcal{O}(X) \rightrightarrows X, U \mapsto U$$

with $\mathrm{dom}(\mathsf{C}_{\mathcal{O}(X)}) := \{U \subseteq X : U \neq \emptyset \text{ open}\}$ is called *open choice* of $X$.

It is clear that open choice is computable, since any non-empty ball $B(x_i, r_i)$ already provides a point in the set, namely its center $x_i$.

COROLLARY 8.2 (Open choice). *Let $X$ be a computable metric space. Then $\mathsf{C}_{\mathcal{O}(X)} :\subseteq \mathcal{O}(X) \rightrightarrows X$ is computable.*

This directly implies the following result, which is nothing but a well-known uniformization property. We include it here in order to emphasize the analogy to the following metatheorems.

THEOREM 8.3 (Open Metatheorem). *Let $X, Y$ be computable metric spaces and let $U \subseteq X \times Y$ be c.e. open. If*

$$(\forall x \in X)(\exists y \in Y)(x, y) \in U,$$

*then $R : X \rightrightarrows Y, x \mapsto \{y \in Y : (x, y) \in U\}$ is computable.*

PROOF. We consider the section map

$$S : X \to \mathcal{O}(Y), x \mapsto U_x := \{y \in Y : (x, y) \in U\}.$$

We use the fact (see [13]) that the representation $\vartheta$ of $\mathcal{O}(X)$ is computably equivalent to the representation $\vartheta'$, defined by

$$\vartheta'(p) = U :\iff [\delta_X \to \rho](p) = f \text{ and } f^{-1}\{0\} = X \setminus U.$$

Since $U \subseteq X \times Y$ is c.e. open, there is a computable function $f : X \times Y \to \mathbb{R}$ such that $(X \times Y) \setminus U = f^{-1}\{0\}$ and by type conversion one obtains that the function

$$g : X \to \mathcal{C}(Y, \mathbb{R}), x \mapsto (y \mapsto f(x, y))$$

is computable as well. However, $(g(x))^{-1}\{0\} = Y \setminus U_x$ and hence $S$ is computable. Thus $R = \mathsf{C}_{\mathcal{O}(Y)} \circ S$ is computable by Corollary 8.2.                    ⊣

An example of a theorem that falls under this category is the Weierstraß Approximation Theorem. It can be formulated as

$$(\forall f \in \mathcal{C}[0, 1])(\forall k \in \mathbb{N})(\exists n \in \mathbb{N}) \ ||f - p_n|| < 2^{-k},$$

where $(p_n)_{n \in \mathbb{N}}$ is some effective enumeration of the rational polynomials $\mathbb{Q}[x]$. The predicate

$$U := \{(f, k, n) : ||f - p_n|| < 2^{-k}\} \subseteq \mathcal{C}[0, 1] \times \mathbb{N} \times \mathbb{N}$$

is c.e. open. Hence, given a continuous function $f : [0, 1] \to \mathbb{R}$ and $k \in \mathbb{N}$ we can actually effectively find a rational polynomial $p_n$ that approximates $f$ with precision $2^{-k}$. If we denote by $\mathsf{WAT} : \mathcal{C}[0, 1] \times \mathbb{N} \rightrightarrows \mathbb{N}$ the corresponding realizer of the Weierstraß Approximation Theorem, then we get the following corollary.

COROLLARY 8.4 (Weierstraß Approximation Theorem). $\mathsf{WAT} \equiv_{\mathrm{W}} \mathrm{id}$.

Of course, this theorem as such is well-known. It has been proved directly, for instance, by Caldwell and Pour-El [16] and Hauck [21]. There are many other approximation results that fall into the same category as the Weierstraß Approximation Theorem. Roughly speaking, if one has a $\Pi_2$–statement with a classical proof such that the corresponding predicate is c.e. open, then one automatically has a computable version of the theorem. The next metatheorem is a similar observation for co-c.e. closed predicates and co-c.e. compact $Y$. We first formulate a lemma.

LEMMA 8.5 (Section). *Let $X, Y$ be computable metric spaces. Then*

$$\mathrm{sec} : \mathcal{A}_-(X \times Y) \times X \to \mathcal{A}_-(Y), (A, x) \mapsto A_x := \{y \in Y : (x, y) \in A\}$$

*is computable.*

The proof is straightforward and can be found in Lemma 6.2 of [8]. As a consequence we obtain the following result.

THEOREM 8.6 (Compact Metatheorem). *Let $X, Y$ be computable metric spaces and let $Y$ be co-c.e. compact and $A \subseteq X \times Y$ co-c.e. closed. If*

$$(\forall x \in X)(\exists y \in Y)(x, y) \in A,$$

*then $R : X \rightrightarrows Y, x \mapsto \{y \in Y : (x, y) \in A\}$ is weakly computable, i.e. $R \leq_W \mathsf{C_K}$.*

PROOF. Because $Y$ is co-c.e. compact, the identity id $: \mathcal{A}_-(Y) \to \mathcal{K}_-(Y)$ is computable (see Lemma 6 in [7]). According to the Section Lemma 8.5 it follows that the map

$$S : X \to \mathcal{K}_-(Y), x \mapsto A_x$$

is computable, because $A \subseteq X \times Y$ is co-c.e. closed. Thus $R = \mathsf{C}_{\mathcal{K}(Y)} \circ S$ is weakly computable by Theorem 2.10.                                    ⊣

In other words, $R \leq_W \mathsf{C_K}$ under the given conditions. An example of a theorem that falls under this category is the Brouwer Fixed Point Theorem which states

$$(\forall f \in \mathcal{C}([0, 1]^n, [0, 1]^n))(\exists x \in [0, 1]^n)\ f(x) = x.$$

The space $[0, 1]^n$ is computably compact and the predicate

$$A = \{(f, x) : f(x) = x\} \subseteq \mathcal{C}([0, 1]^n, [0, 1]^n)) \times [0, 1]^n$$

is co-c.e. closed. Hence there is a weakly computable solution for the Brouwer Fixed Point Theorem. If we denote by $\mathsf{BFT} : \mathcal{C}([0, 1]^n, [0, 1]^n)) \rightrightarrows [0, 1]^n$ the corresponding realizer of the Brouwer Fixed Point Theorem, then we get the following corollary.

COROLLARY 8.7 (Brouwer Fixed Point Theorem). $\mathsf{BFT} \leq_W \mathsf{WKL}$.

Having such an upper bound allows immediately to draw a lot of conclusions. From the Low Basis Theorem of Jockusch and Soare and the Low Invariance Principle 4.4 we get for instance the following conclusion.

COROLLARY 8.8. *Let $f : [0, 1]^n \to [0, 1]^n$ be a computable function. Then $f$ has some low fixed point $x \in [0, 1]^n$, i.e. $f(x) = x$ and $x$ has some low name $p \in \mathbb{N}^{\mathbb{N}}$ such that $\rho(p) = x$.*

Many other theorems of analysis that have to do with the solution of equations in compact spaces fall into the same category. This applies for instance to the Schauder Fixed Point Theorem and also to the Intermediate Value Theorem. Sometimes it is not immediately clear that a theorem is of this form. In case of the Peano Existence Theorem for solutions of initial value problems of ordinary differential equations it is easy to see that it can be reduced to the Schauder Fixed Point Theorem (see [48]). Another example of this type is the Hahn-Banach Theorem. As it is usually formulated, is not of the form of an equation with a solution in a compact space. However, using the Banach-Alaoglu Theorem, it can be brought into this form (see [7, 20]). Whenever a theorem that falls under the Compact Metatheorem has a unique solution, then that solution is automatically uniformly computable in the input. That follows from the Compact Metatheorem 8.6 and Corollary 2.11, which states that any single-valued weakly computable function is computable.

COROLLARY 8.9 (Unique Compact Metatheorem). *Let $X, Y$ be computable metric spaces. Let $Y$ be co-c.e. compact and let $A \subseteq X \times Y$ be co-c.e. closed. If*

$$(\forall x \in X)(\exists! y \in Y)(x, y) \in A,$$

*then $R : X \to Y, x \mapsto \{y \in Y : (x, y) \in A\}$ is computable.*

Thus, under all (perhaps purely classical) conditions under which the Brouwer Fixed Point Theorem, the Intermediate Valued Theorem, the Hahn-Banach Theorem or the Peano Existence Theorem have unique solutions, they are already automatically fully computable. Finally, we want to prove a locally compact version of our metatheorem. As a preparation we first define a generalization of closed choice.

DEFINITION 8.10 (Closed choice). Let $X$ be a computable metric space. The multi-valued operation

$$\mathsf{C}_{\mathcal{A}(X)} :\subseteq \mathcal{A}_-(X) \rightrightarrows X, A \mapsto A$$

with $\mathrm{dom}(\mathsf{C}_{\mathcal{A}(X)}) := \{A \subseteq X : A \neq \emptyset \text{ closed}\}$ is called *closed choice* of $X$.

In general closed choice is much less well-behaved than compact choice. It is known, for instance by a result of Kleene, that there are co-c.e. closed sets $A \subseteq \mathbb{N}^{\mathbb{N}}$ that have no hyperarithmetical points (see [29] and Theorem 1.7.1 in [17]). A relativization of that result leads to the following conclusion.

COROLLARY 8.11. $\mathsf{C}_{\mathcal{A}(\mathbb{N}^{\mathbb{N}})} :\subseteq \mathcal{A}_-(\mathbb{N}^{\mathbb{N}}) \rightrightarrows \mathbb{N}^{\mathbb{N}}$ *has no Borel measurable realizer.*

On the other hand, we know that $\mathsf{C}_{\mathcal{A}(\{0,1\}^{\mathbb{N}})} \equiv_{\mathrm{W}} \mathsf{C}_{\mathsf{K}} <_{\mathrm{W}} \mathsf{C}_{\mathsf{A}} = \mathsf{C}_{\mathcal{A}(\mathbb{R})}$ by Corollary 4.12 and Theorem 2.10 and $\mathsf{C}_{\mathcal{A}(\mathbb{R})}$ is $\mathbf{\Sigma}_2^0$–computable (but not $\mathbf{\Sigma}_2^0$–complete). Thus, the Weihrauch degree of $\mathsf{C}_{\mathcal{A}(X)}$ sensitively depends on $X$. However, if the computable metric space $X$ is effectively locally compact, then we can say at least something. We recall that $X$ is *effectively locally compact*, if there is an operation that computes for any point $x \in X$ and any rational open neighbourhood $I$ of $x$ some compact set $K$ with full $\kappa$–information such that $x \in K^\circ \subseteq K \subseteq I$ (see [10] for details). In Proposition 4.5 of [10] it has been proved that for effectively locally compact $X$ the identity $\mathrm{id} : \mathcal{A}_-(X) \to \mathcal{A}_+(X)$ is limit computable (i.e. $\mathbf{\Sigma}_2^0$–computable). Since any effectively locally compact metric space $X$ is complete, positive information on non-empty closed sets allows to select a point (see [13]). We obtain the following corollary.

COROLLARY 8.12 (Closed Choice). *Let $X$ be an effectively locally compact computable metric space. Then $\mathsf{C}_{\mathcal{A}(X)} \leq_{\mathrm{W}} \mathsf{C}$.*

This corollary allows us to prove the following locally compact version of our Metatheorem.

THEOREM 8.13 (Locally Compact Metatheorem). *Let $X, Y$ be computable metric spaces, let $Y$ be effectively locally compact and let $A \subseteq X \times Y$ be co-c.e. closed. If*

$$(\forall x \in X)(\exists y \in Y)(x, y) \in A,$$

*then $R : X \rightrightarrows Y, x \mapsto \{y \in Y : (x, y) \in A\}$ satisfies $R \leq_{\mathrm{W}} \mathsf{C}_{\mathcal{A}(Y)}$. In particular, $R$ is limit computable.*

Proof. According to the Section Lemma 8.5 the map

$$S : X \to \mathcal{A}_-(Y), x \mapsto A_x$$

is computable with respect to $\psi_-$, because $A \subseteq X \times Y$ is co-c.e. closed. According to Corollary 8.12 the closed choice $\mathsf{C}_{\mathcal{A}(Y)}$ is limit computable because $Y$ is effectively locally compact. Thus $R = \mathsf{C}_{\mathcal{A}(Y)} \circ S \leq_{\mathrm{W}} \mathsf{C}_{\mathcal{A}(Y)}$ and $R$ is limit computable too.                                                                    ⊣

We illustrate an application of the Locally Compact Metatheorem.

Corollary 8.14. *Let $X$ be a computable metric space and let $Y$ be an effectively locally compact metric space. If $f : X \to Y$ is a function with a co-c.e. closed graph* $\mathrm{graph}(f) = \{(x, y) \in X \times Y : f(x) = y\}$, *then $f$ is limit computable. In particular, the inverse $g^{-1} : X \to Y$ of any computable bijective function $g : Y \to X$ is limit computable.*

One can even say more in this case. The function $f|_K$ restricted to any compact $K \subseteq X$ is computable (see Proposition 7.2 in [8]) and hence $f$ maps computable points to computable points. We now get a corresponding uniqueness version of the Locally Compact Metatheorem as a corollary.

Corollary 8.15 (Unique Locally Compact Metatheorem). *Let $X, Y$ be computable metric spaces, let $Y$ be effectively locally compact and let $A \subseteq X \times Y$ be co-c.e. closed. If*

$$(\forall x \in X)(\exists! y \in Y)(x, y) \in A,$$

*then $R : X \to Y, x \mapsto \{y \in Y : (x, y) \in A\}$ satisfies $R \leq_{\mathrm{W}} \mathsf{C}_{\mathcal{A}(Y)}$. In particular, $R$ is limit computable and maps computable inputs to computable outputs. The function $R|_K$ restricted to compact sets $K \subseteq X$ is computable.*

§9. **Conclusions.** In this paper we have suggested a new approach to classify mathematical theorems according to their computational content just by using methods of topology and computability theory. It is obvious that this approach is closely related to other classifications of more proof theoretic nature that exist in constructive and reverse mathematics. We cannot provide any exhaustive analysis of the similarities and differences between our results and known results, but we briefly summarize some observations.

**9.1. Constructive mathematics.** Based on the work of Bishop and Bridges [1, 15] Ishihara has classified many theorems in constructive analysis with respect to their relation to certain non-constructive principles, in particular, with respect to LPO and LLPO, see the survey [25]. Equivalence of theorems in this approach essentially means intuitionistical equivalence in a setting where certain choice axioms such as countable choice, dependent choice and unique choice are accepted. It turns out that in this setting the following theorems are equivalent to LLPO:

1. Weak Kőnig's Lemma [24].
2. The Intermediate Value Theorem [15].
3. The Hahn-Banach Theorem [24].
4. Some form of the Heine-Borel Theorem [24].

One of the several differences to our classifications is that the Heine-Borel Theorem is fully computable in our approach [52]. This is essentially because we do not have to prove existence of finite subcovers, but we can just search for them systematically (which is also impossible in constructive analysis without Markov's principle).

Moreover, it is clear that there is no distinction between LLPO and its parallelization $\widehat{\text{LLPO}}$ in constructive analysis, since the parallel and even the sequential application of $\widehat{\text{LLPO}}$ is allowed in the logical framework of constructive analysis. The sequential application of $\widehat{\text{LLPO}}$ does not lead to any more difficult operations, since the class is closed under composition (see [20] and [12]).

This is different for LPO, since $\widehat{\text{LPO}}$ is known to be $\mathbf{\Sigma}_2^0$–complete (see [5]) it follows that the iterations of this operation climb up the finite part of the Borel hierarchy. This is why equivalence to LPO in the framework of constructive analysis is a much wider concept and even includes the Bolzano-Weierstraß Theorem [25, 35], which is not below $\widehat{\text{LPO}}$ in our framework (this will be shown elsewhere).

Other differences between constructive analysis and our approach are that we cannot distinguish between LPO and its weaker version WLPO, since both principles are equivalent in our approach and indeed their realizations have the same degree of discontinuity. Principles like Markov's principle MP are continuous and computable from our perspective.

Other weaker principles of omniscience such as $\text{LLPO}_n$ that have been introduced in constructive analysis by Richman [41, 42] have, however, been intensively studied in computable analysis [51, 36] and they seem to be useful for the classification of problems of more combinatorial nature.

Finally, we mention that our uniqueness results for weakly computable functions are reminiscent of the study of uniqueness questions in constructive analysis by Schuster, see [45, 46] and our metatheorems are perhaps related to those of Gelfond, see [18].

**9.2. Reverse mathematics.** In reverse mathematics as proposed by Friedman and Simpson [48] theorems are classified according to which comprehension axioms are required to prove the corresponding theorems in second order arithmetic. Almost all theorems that we have considered have also been classified in reverse mathematics. For instance the following theorems are known to be provable in the base system $\text{RCA}_0$ with recursive comprehension (see [48]):

1. The Baire Category Theorem.
2. The Intermediate Value Theorem.
3. The Uniform Boundedness Theorem.

There are many theorems provable in the base system that are uniformly computable in our approach, such as the Tietze Extension Theorem [53] and we do not list them there. In case of the Baire Category Theorem the actual proof provided in [48] is a proof of the constructive and computable version $\text{BCT}_0$. But since reverse mathematics is based on classical logic this statement can be freely converted into the non-computable and non-constructive version BCT. Correspondingly, all theorems that are derivable from BCT are also provable in $\text{RCA}_0$. Roughly speaking, the base system $\text{RCA}_0$ corresponds to our non-uniform computability results.

The second most important system in reverse mathematics is $\mathsf{WKL}_0$ and it roughly corresponds to our class of weakly computable operations (operations that are Weihrauch reducible to $\mathsf{WKL}$). This class includes the Hahn-Banach Theorem and, of course, Weak Kőnig's Lemma itself [48]. Similarly to constructive mathematics (see above) we find that the Heine-Borel Theorem is a distinguishing feature, as it is computable in our system. Finally, the system $\mathsf{ACA}_0$ of arithmetic comprehension corresponds to our class of effectively Borel measurable maps and in contrast to our system composition of operations is always for free like in constructive mathematics which makes it difficult to distinguish different levels of the Borel hierarchy in reverse mathematics. In contrast to constructive mathematics, however, parallelization does not seem to be for free in reverse mathematics. The relation between computable analysis and reverse mathematics has first been studied in [20].

Kohlenbach has suggested a uniform version of reverse mathematics [32, 33, 44]. Although the uniform formulations of the studied principles such as Weak Kőnig's Lemma seem to be very closely related to our realizer interpretations, the results are considerably different. For instance, it seems that the analogue $(\exists^2)$ of $\mathsf{LPO}$ is equivalent to the uniform version of Weak Kőnig's Lemma, whereas in our approach these two principles are even incomparable.

Moreover, we mention that there could be a relation between known conservativeness results of $\mathsf{WKL}_0$ in reverse mathematics [47] and the fact that single-valued weakly computable function are computable in our setting.

Finally, we note that there is the approach of limit computable mathematics proposed by Hayashi and studied by several others. Also in this context a close relation between $\mathsf{WKL}$ and $\mathsf{LLPO}$ has been established [37, 22].

The observations collected here lead to plenty of fascinating questions for further research. In any case, we believe that our computational classification of theorems contributes new insights and yields a finer classification of theorems in some respects. To understand the exact relation between our approach, the intuitionistic approach in constructive analysis and the more proof theoretic analysis in reverse mathematics requires further studies.

LABORATORY OF FOUNDATIONAL ASPECTS OF COMPUTER SCIENCE
    DEPARTMENT OF MATHEMATICS & APPLIED MATHEMATICS
        UNIVERSITY OF CAPE TOWN
            RONDEBOSCH 7701, SOUTH AFRICA

DIPARTIMENTO DI FILOSOFIA
    UNIVERSITÀ DI BOLOGNA
        ITALY
*E-mail*: Vasco.Brattka@uct.ac.za
*E-mail*: Guido.Gherardi@unibo.it
*URL*: http://cca-net.de/